\documentclass{amsart}
\usepackage{amsfonts}
\usepackage{amssymb}
\begin{document}

{\large\bf
Congruence conditions, parcels, and Tutte polynomials of graphs and matroids}

\vskip .33in
\ \ Joseph P.S. Kung\footnote{ Supported by the National Security Agency under Grant H98230-11-1-0183.}

\ \ Department of Mathematics,

\ \ University of North Texas, Denton, TX 76203, U.S.A.

\vskip .33in 
\ \ e-mail: kung@unt.edu
\vskip .50in

\subjclass{Primary 05B35;
Secondary 05B20, 05C35, 05D99, 06C10, 51M04}

\newcommand\join{\vee}
\newcommand\meet{\wedge}
\newcommand\rk{\mathrm{rk}}
\newcommand\aA{\mathbb{A}}
\newcommand\wt{\mathrm{wt}}
\newcommand\RR{\mathrm{R}}
\newcommand\ee{\mathrm{e}}

\noindent
{\it ABSTRACT.} 

Let $G$ be a matrix and $M(G)$ be the matroid defined by linear dependence on the set $E$ of column vectors of $G.$  Roughly speaking, a parcel is a subset of pairs 
$(f,g)$ of functions defined on $E$ to a suitable Abelian group $\mathbb{A}$ satisfying a coboundary condition (that 
the difference $f-g$ is a flow over $\mathbb{A}$ of $G$) and a congruence condition (that an algebraic or combinatorial function of $f$ and $g,$ such as the sum of the size of the supports of $f$ and $g,$ 
satisfies some congruence condition).    
We prove several theorems of the form:  a linear combination of sizes of parcels, with coefficients roots of 
unity, equals a multiple of an evaluation of the Tutte polynomial of $M(G)$ at a point $(u,v),$
usually with complex coordinates, satisfying $(u-1)(v-1) = |\mathbb{A}|.$    

\vskip .33in \noindent
{\it Subject classification.}  05B35  05C15  05C80

\noindent
{\it Keywords.}  Flow, matroid, parcel, Tutte polynomial

\newpage
\vskip 1.2in \noindent
{\large \bf 1.  Flows in graphs and matrices} 

\vskip .20in \noindent
There are several results in graph theory saying {\it in essence} that a multiple of the chromatic, flow, or Tutte polynomial evaluated 
at a certain point equals the difference in size of two sets or ``parcels'' of pairs or triples of subsets of the edge set 
satisfying a coboundary condition and different congruence conditions on their size \cite{Goodall1, Goodall2, Goodall, Mat}.   
A simple example of such a result is the following analog of Theorem 18 in Goodall's paper \cite{Goodall}. (This analog will be proved as a special case of Corollary 4.15 in this paper.)      

\vskip 0.2in\noindent 
{\bf Theorem 1.1.}
Let $\Gamma (V,E)$ be a graph on the vertex set $V$ and the edge set $E$ with $c$ connected components.  
 For $k = 0,1,2,3,$ let 
${\mathcal{L}}(k)$ (respectively, $\bar{\mathcal{L}}(k)$) be the number of pairs $(A,B)$ of edge-subsets such that 
the subgraph with edge set the symmetric difference $A \Delta B$  
is a disjoint union of minimal cutsets (respectively, is a disjoint union of cycles) and $|A| + |B| \equiv |E| + k$ modulo $ 4.$    Then 
$$
|{\mathcal{L}}(1)| = |{\mathcal{L}}(3)| \,\,\,\,\mathrm{and}\,\,\,\, 
|{\mathcal{L}}(0)|  - |{\mathcal{L}}(2)|  = 2^{|E|-c} P(\Gamma;2), 
$$
where $P(\Gamma;\lambda)$ is the chromatic polynomial of $\Gamma,$ and   
$$
|\bar{\mathcal{L}}(1)| = |\bar{\mathcal{L}}(3)| \,\,\,\,\mathrm{and}\,\,\,\, 
|\bar{\mathcal{L}}(0)|  - |\bar{\mathcal{L}}(2)|  = 2^{|E|} F(\Gamma;2), 
$$
where $F(\Gamma;\lambda)$ is the flow polynomial of $\Gamma.$    

\vskip 0.3in
An aim of the graph-theoretic results is to reformulate theorems, like the $4$-color theorem, in probabilistic terms.  
The statements of these results involve differences in probabilities and their proofs are intricate, involving commutative algebra or finite Fourier analysis.  Our aim in this paper is different: it is to put these results in context, as initial cases of infinite families of results about sizes of parcels constructed using flows of matrices over suitable Abelian groups.    

This paper is the second in a series.  The earlier paper \cite{cob} is about the ``parametric'' theory.  Using characters, 
we studied parcels defined by algebraic conditions using actual values of elements in the Abelian group $\aA.$  
With two exceptions, the present paper is about the ``non-parametric'' theory.  Parcels are defined using congruence conditions which only use the data whether an element in $\aA$ is zero or nonzero.   We will use non-parametric parcels to obtain combinatorial interpretations of evaluations of rank generating polynomials on complex hyperbolas $\lambda x = q,$   
where $q$ is an  integer greater than $1.$     
(There seems to be only one interpretation of the evaluation of the rank generating polynomial at complex points in the  literature, the interpretation by Jaeger \cite{Jg} of the value of the Tutte polynomial at $(\ee^{2\pi \iota/3},\ee^{-2 \pi \iota/3}).$ Jaeger's interpretation will be discussed in Section 6.)   

We shall assume some knowledge of graph and matroid theory.  See, for example, \cite{BryOx, KungC, Oxley}.
We recall briefly concepts and definitions central to this paper.   
Let $G$ be a matrix with columns indexed by the set $E$ and $M(G)$ be the matroid on $E$ defined by linear dependence of the column vectors.  Let $\Gamma (V,E)$ be a graph with a chosen orientation on the 
edges.  The orientation gives a function $\mathrm{sign}(v,e):$  if $e$ is not a loop, then $\mathrm{sign}(v,e)$ equals 
$+1,$ $-1,$ or $0,$ depending on whether 
the edge $e$ is going into $v,$ going out of $v,$ or not incident at all on $v.$  If $e$ is a loop, then $\mathrm{sign}(v,e)$ 
is always $0.$    
We will use two matrices defined by $\Gamma.$  The first is the {\sl vertex-edge matrix}, the 
$|V| \times |E|$ matrix $H$ with rows indexed by $V$ and columns indexed by $E$ such that the $ve$-entry 
is $\mathrm{sign}(v,e).$   The matroid on $E$ defined by the vertex-edge matrix is the cycle matroid of $\Gamma.$  The second is the cycle-edge matrix.  We think of a cycle $c$ as a set of edges 
$e_0,e_1, \ldots, e_{t-1}$ such that there is a sequence 
$v_0,v_1, \ldots,v_{t-1}$ of distinct vertices with the edge $e_j$ incident on $v_j$ and $v_{j+1}.$   The indices are 
regarded as integers modulo $t.$  The {\sl cycle-edge matrix} $G$ has rows indexed by the set of all cycles, columns indexed by $E,$ and $ce$-entry equal to $0$ unless $e$ equals one of the edges $e_j$ in the cycle $c,$ in which case it equals $\mathrm{sign}(v_j,e_j).$   The matroid on $E$ defined by the cycle-edge matrix is 
the cocycle matroid of $\Gamma.$   

Two row vectors  $(u_e)$ and $(v_e)$ (indexed by the set $E$) are {\sl orthogonal} if their inner product 
$\sum_{e:\,e \in E} u_e v_e$ equals zero.   
Two matrices $H$ and $G$ on the same column set $E$ are {\sl orthogonal duals} of each other if every row of $H$ is orthogonal to every row of $G$ and $\mathrm{rank}(H) + \mathrm{rank}(G) = |E|.$   
For example, the cycle-edge matrix $G$ and the vertex-edge matrix $H$ of a graph $\Gamma$ are orthogonal duals. 
The matroid $M(H)$ is the matroid  $M(G)^{\perp}$ dual to $M(G).$  

Let $\aA$ be an additive Abelian group.  A {\sl flow} $h$ ({\sl over $\aA$}) on the graph $\Gamma$ is a function $h:E \to \aA$ such that at each vertex, the inflow equals the outflow.  To state this {\sl conservation condition} precisely, associate with the function $h:E \to \aA$  the {\sl row vector} $(h(e))_{e \in E}$ with coordinates in $\aA$ indexed by $E$ and denote both the function and its row vector by $h.$  Then the conservation condition at all the vertices can be stated succinctly by the matrix equation  
$$
H h^T = 0, 
$$ 
where $0$ is the zero column vector, and $\cdot^T$ is transpose, so that $h^T$ is a column vector.  Since $H$ is a matrix with entries in $\{-1,0,1\},$ the matrix product is defined. 

Flows over fields can be defined over any pairs of orthogonal duals.  Let $G$ be a rank-$r$ $n \times |E|$ matrix 
with columns indexed by $E$ with entries in a field $\mathbb{F}.$  Then by solving a system of linear equations, we 
can find a matrix $H$ such that $G$ and $H$ are orthogonal duals.  
We define a {\sl flow} $h$  of $G$ {\sl over} $\mathbb{F}$   
to be a function $h: \,E \to \mathbb{F}$ such that $H h^T = 0.$   From linear algebra, a row vector $h$ satisfies $Hh^T =0$ if and only if it is a linear combination of row vectors of $G.$   
Hence, as a row vector, $h$ is a flow over $G$ if and only if it is in the row space of $G.$   
Note that in the case when $G$ is a $0 \times |E|$ matrix (and the matroid $M(G)$ has rank $0$), there is a unique flow over $G,$ the zero row vector.      

A {\sl linear functional} is a linear function $\mathbb{F}^n \to \mathbb{F}.$   If $e$ is in the column set of an 
$n \times |E|$ matrix $G,$ then $\ell(e)$ is the value of $\ell$ on the column vector indexed by $e.$     

\vskip 0.2in \noindent 
{\bf  Lemma 1.2.} Let $h$ be a flow of the matrix $G$ over $\mathbb{F}.$  Then there exists a linear 
functional $\ell$ such that for all $e$ in $E,$  $h(e) = \ell(e).$
The linear functional $\ell$ is uniquely determined by $h$ if $n=r.$   
\vskip 0.1in \noindent 
{\bf  Proof.}   Let  
$g_1,g_2,\ldots,g_n$ be the rows of $G$ and $h = a_1 g_1 + a_2 g_2+ \cdots + a_n g_n.$
Then the linear functional $\ell_h$ sending the column vector $(x_1,x_2,\ldots,x_n)^T$ to 
$a_1x_1 + a_2x_2 + \cdots + a_nx_n$ has the required property that for all $e,$ $\ell_h(e)= h(e).$  
\qed

\vskip 0.2in
The description of a flow as a linear combination of row vectors carries over to graphs.  A flow $h$ on a graph $\Gamma$ over an Abelian group $\aA$ is {\sl minimal} if there is a nonzero group element $a$ and a cycle $c$ of $\Gamma$ such that $h(e) = 0$ if $e$ is not in $c$ and $h(e_j) = \mathrm{sign}(v_j,e_j)a$ for the edges $e_j$ in $c.$   It is easy to prove that $h$ is a flow over a graph if and only if it is a finite sum of minimal flows.  Put another way, a flow is a ``linear combination'' of row vectors of the cycle-edge matrix with coefficients in $\aA.$   
One can dualize the notion of flows on graphs and define a {\sl tension} $h^{\prime}$ on a graph $\Gamma$ to be a function $h^{\prime}: E \to \aA$ such that 
$G h^{\prime T} = 0,$ or equivalently, the row vector $h^{\prime}$ a ``linear combination'' of the rows of the vertex-edge matrix with coefficients in $\aA.$ 
A uniform treatment of flows and tensions can be obtained by viewing them as flows on an orthogonally dual  pair of totally unimodular matrices.      

Let $f: E \to \aA$ be a function.  The {\sl kernel} of $f$ is the inverse image $f^{-1}(0)$ in $E.$  

\vskip 0.2in \noindent 
{\bf  Lemma 1.3.} 
The kernel of a flow of the matrix $G$ is 
a closed set in the matroid $M(G).$   
\vskip 0.2in \noindent 
{\bf  Proof.} 
For a flow $h$ over a field, $h^{-1}(0)$ is the intersection of $E$ with the null space of the linear function $\ell_h.$  Hence, $h^{-1}(0)$ is closed under linear dependence.   For a flow on a graph, the kernel is the complement of an edge-disjoint union   
of cycles and hence, a closed set in the cocycle matroid.  For a tension on a graph, the kernel is the complement of an edge-disjoint union of minimal cutsets and hence, a closed set in the cycle matroid. 
\qed
\vskip 0.2in 

The {\sl support} $\mathrm{supp}(f)$ of a function $f: E \to \aA$ is the complement $E \backslash h^{-1}(0);$ in other words,  
$$
\mathrm{supp}(f) = \{e: \, f(e) \neq 0\}.
$$
The support of a flow is the complement of a closed set and hence, it is a union of cocircuits.  The minimal subsets among 
supports of flows are exactly the cocircuits of the matroid $M(G).$  
In the case $h$ is a flow over a field $\mathbb{F},$ $\mathrm{supp}(h)$ is contained in the complement of the hyperplane in $\mathbb{F}^n$ on which the linear functional  $\ell_h$ is zero; in other words, $\mathrm{supp}(h)$ is {\sl affine}
over $\mathbb{F}.$

The field $\mathrm{GF}(2)$ of order $2$ has exactly one nonzero element.  Hence, a function $f:E \to \mathrm{GF}(2)$  is determined by its support and one can identify $f$ with the subset $\mathrm{supp}(f)$ in $E.$    For easy reference, we recall the following elementary results (see, for example, Propositions 9.3.1 and 9.3.2 in \cite{Oxley}).  

\vskip 0.2in \noindent 
{\bf  Lemma 1.4.}    
Let $G$ be a binary matrix with column set $E$ and $B \subseteq E.$  
The following are equivalent:  the subset $B$ is $\mathrm{GF}(2)$-affine; all circuits contained in 
$B$ have even size;  $\chi(M|B;2) = 1,$ where $M|B$ is the restriction of $M$ to $B.$  

\vskip 0.2in \noindent 
{\bf  Lemma 1.5.} 
Let $\Gamma(V,E)$ be a graph and $h:E \to \mathrm{GF}(2)$ be a function.  Then 

(a) $h$ is a flow if and only if the subgraph on the edge-subset $\mathrm{supp}(h)$ is a (disjoint) union of cycles, or equivalently, all the vertices on that subgraph have even degree.

(b)  $h$ is a tension if and only if the subgraph on $\mathrm{supp}(h)$ is a (disjoint) union of minimal cutsets.   

\vskip 0.2in\noindent
The two cases, graphs and matrices, call for a theory encompassing both.    
A reasonable theory that covers all known cases is the theory of totally-$\mathbb{P}$ matrices and $\mathbb{P}$-modules, 
where $\mathbb{P}$ is a (concrete) partial field \cite{Semple}.   This theory is described in \cite{cob} and we shall not repeat the description here.   However, we shall use the technical hypothesis 
``Let $G$ be a totally $\mathbb{P}$-matrix and $\aA$ be a 
$\mathbb{P}$-module of order $q.$''   
The reader unfamiliar with (or unconvinced by) totally-$\mathbb{P}$ matrices should read the hypothesis 
as a short way to write ``Let $G$ be the vertex-edge matrix or  
the cycle-edge matrix of a graph, or a totally unimodular matrix, and $\aA$ be an Abelian group of order $q,$ 
or let $G$ be a matrix over $\mathrm{GF}(p^s),$ $\aA$ be the vector space $\mathrm{GF}(p^s)^d,$ and   
$q = p^{sd}.$'' 

We begin Section 2 by explaining how two identities due to Tutte can be used to obtain a  third identity relating weights of $(m+1)$-tuples of functions to a sum over flows.  Parcels are defined next.  In many cases, the third identity specializes to a relation between sizes of parcels and an evaluation of the rank generating polynomial.  We then apply this theory to parcels defined on pairs of functions by Hamming distances (Section 3), sizes of supports (Section 4), and inner products (Section 6).  Parcels involving triples and $(m+1)$-tuples of functions are discussed in Section 5.  Finally, in Section 7, we discuss generic or enumerator versions of our results.  

We shall work in the field $\mathbb{C}$ of complex numbers with the notation: $\ee = 2.7182\ldots$ and 
$\iota = \sqrt{-1}.$

\vskip .5in \noindent 
{\bf\large  2.  The Cheshire-cat identity}

\vskip .20in \noindent
Let $M$ be a rank-$r$ matroid on the set $E$ with rank function $\rk.$  
The {\sl rank generating} or {\sl corank-nullity polynomial} of $M$ is the polynomial $R(M;\lambda,x)$ in the variables $\lambda$ and $x$ 
defined by 
$$
R(M;\lambda,x) = \sum_{B:\, B \subseteq E}  \lambda^{r - \rk(B)} x^{|B| - \rk (B)}.  
$$     
The rank generating polynomial satisfies the duality condition:  if $M^{\perp}$ is the dual of $M,$ then 
$ R(M^{\perp} ;\lambda,x) = R(M;x,\lambda). $
The {\sl characteristic polynomial} $\chi (M;\lambda)$ of $M$ is defined to be the polynomial 
$$
\sum_{B: B \subseteq E} (-1)^{|B|}  \lambda^{r - \rk(B)}.  
$$
Up to a sign, the characteristic polynomial is an evaluation of the rank generating polynomial: indeed, 
$$
\chi (M;\lambda) = (-1)^{r} R(M;-\lambda,-1) = (-1)^{r} R(M^{\perp};-1,-\lambda) .  
\eqno(1)$$
The {\sl Tutte polynomial} $T(M;u,v)$ is defined to be the polynomial $R(M;u-1,v-1).$  Despite the title, we will use the rank generating polynomial instead of the Tutte polynomial.   We shall need the following identity of Crapo and Tutte \cite{Crapo, Tutte}.  

\vskip 0.3in\noindent 
{\bf Lemma 2.1.}   Let $M$ be a rank-$r$ matroid with lattice $L(M)$ of flats.  Then 
$$ (x-1)^r R\left(M; \frac {\lambda}{x-1},x-1\right) = \sum_{U:\, U\in L(M)} \chi(M/U; \lambda) x^{|U|}. $$ 

\vskip 0.2in 
When $H$ is the vertex-edge matrix of a graph $\Gamma$ with $c$ connected components, then the chromatic polynomial 
$P(\Gamma;\lambda)$ equals $\lambda^c\chi(M(H);\lambda).$  When $G$ is the cycle-edge matrix 
of $\Gamma,$ then the flow polynomial $F(\Gamma ; \lambda)$ equals $\chi(M(G);\lambda).$ 
The critical problem (\cite{CR}; see also \cite{KungC}) of Crapo and Rota gives a counting interpretation 
of the characteristic polynomial extending the interpretations for graphs.  We  state this interpretation in the following general form \cite{cob}.  Let $G$ be a totally-$\mathbb{P}$ matrix with column set $E,$ 
$\aA$ be a $\mathbb{P}$-module of order $q,$ and $U \subseteq E.$  Then  
$$
\chi (M(G)/U; q) = |\{h:\, h \,\,\text{is a flow over}\,\,\aA, \, h^{-1}(0) = U \}|, 
$$
where $M(G)/U$ is the contraction of the matroid $M(G)$ by $U.$  In particular, $\chi(M(G);q)$ is the number of 
nowhere-zero flows over $\aA.$  Together with Lemma 2.1, this interpretation yields another identity of Tutte \cite{Tutte} as generalized in \cite{cob}.    

\vskip 0.2in\noindent 
{\bf Lemma 2.2.}   Let $G$ be a rank-$r$ totally-$\mathbb{P}$ matrix and $\aA$ a $\mathbb{P}$-module of order $q.$  Then 
\begin{eqnarray*}
\sum_{h:\,h \, \mathrm{is}\,\, \mathrm{a}\,\, \mathrm{flow}\,\,\mathrm{over}\,\,\aA}  x^{|h^{-1}(0)|} 
&=& \sum_{U:\, U \in L(M)}  \chi (M(G)/U; q) x^{|U|} 
\\
&=&  (x-1)^r R\left(M(G);\frac {q}{x-1},x-1\right).
\end{eqnarray*}
  
\vskip 0.3in
We can now state and prove the {\sl Cheshire-cat identity.}   This identity generalizes identities in \cite{Goodall1, cob}.

\vskip 0.2in\noindent 
{\bf Lemma 2.3.}
Let $G$ be a totally-$\mathbb{P}$ matrix, $\aA$ be a $\mathbb{P}$-module, and $\mathbb{B} \subseteq \aA.$  For each 
$(m+1)$-tuple $(b_1,b_2, \ldots,b_{m+1})$ in $\mathbb{B}^{m+1},$ 
let $\gamma (b_1,b_2, \ldots,b_{m+1})$ be an indeterminate.  Then 
\begin{eqnarray*} 
&& \sum_{(f_1,f_2,\ldots,f_{m+1})}
%
%
\quad\quad  \prod_{e:\, e \in E}  \gamma(f_1(e),f_2(e),\ldots,f_{m+1}(e))
\\
&=& 
\sum_{(h_1,h_2,\ldots,h_m)}
\quad
\prod_{e:\,e \in E} \,\,\left(
\sum_
{(b_1,b_2,\ldots,b_{m+1}): \, b_j \in \mathbb{B},\, b_j - b_{j+1} = h_j(e)}
\gamma (b_1,b_2,\ldots,b_{m+1})
\right),
\end{eqnarray*} 
where the sum on the left-hand side ranges over all $(m+1)$-tuples 
$(f_1,f_2,\ldots,f_{m+1})$ of functions $f_j:\, E \to \mathbb{B}$ such that $f_1 - f_2, f_2 - f_ 3, \ldots, f_m - f_{m+1}$ are flows, the outer sum on the right-hand side ranges over all $m$-tuple $(h_1,h_2,\ldots,h_m)$  of flows over $\aA,$ and 
the inner sum on the right-hand side ranges over all $(m+1)$-tuples $(b_1,b_2,\ldots,b_{m+1})$ such that for $1 \leq j \leq m,$ $b_j \in \mathbb{B}$ and $b_j - b_{j+1} = h_j(e).$     

\vskip 0.2in\noindent 
We remark that since the sum of two flows is a flow, the condition that the differences $f_1 -f_2, f_2 -f_3, \ldots, f_m -f_{m+1}$ are flows is equivalent to the condition that $f_j - f_l$ are flows for all indices $j$ and $l.$  Both conditions are equivalent to the condition that the $m$-tuple $( f_1 -f_2, f_2 -f_3, \ldots, f_m -f_{m+1})$ is a flow over $\aA^m.$ 

\vskip 0.2in\noindent 
{\bf Proof of Lemma 2.3.}
Let $(h_1,h_2,\ldots,h_m)$ be a given $m$-tuple of flows.  Expanding the product 
$$ 
\prod_{e:\,e \in E} \,\,\left(
\sum_
{(b_1,b_2,\ldots,b_{m+1}):\,b_j \in \mathbb{B},\, b_j - b_{j+1} = h_j(e)}
\gamma (b_1,b_2,\ldots,b_{m+1}) \right).  
$$
on the right-hand side by the distributive law, we obtain a sum of monomials.  The monomials are products   
$$
\prod_{e:\, e \in E} \gamma (b_1^e,b_2^e,\ldots,b_{m+1}^e)      
$$
with one indeterminate $\gamma (b_1^{e},b_2^{e},\ldots,b_{m+1}^{e})$ for each edge $e,$ and for that indeterminate, $b_j^e \in \mathbb{B}$ and $b_j^e - b_{j+1}^e = h_j(e).$  With each monomial, we associate the unique $(m+1)$-tuple of functions $(f_1,f_2,\ldots,f_{m+1})$ defined by 
$$
f_j (e) =b_j^e. 
$$
By definition, $f_j:\,E \to \mathbb{B}$ and $f_j - f_{j+1} = h_j.$  Conversely, each $(m+1)$-tuple $(f_1,f_2,\ldots,f_{m+1})$ of functions $f_j:\,E \to \mathbb{B}$ such that $f_j - f_{j+1} = h_j$ contributes the monomial
$$
\prod_{e:\,e\in E}  \gamma (f_1 (e), f_2(e), \ldots,f_{m+1}(e)).   
$$
Hence,   
\begin{align*}  
& \prod_{e:\,e \in E} \,\,\left(
\sum_
{(b_1,b_2,\ldots,b_{m+1}):\,b_j \in \mathbb{B},\, b_j - b_{j+1} = h_j(e)}
\gamma (b_1,b_2,\ldots,b_{m+1}) \right)  
\\
& \qquad\qquad  
= \sum_{(f_1,f_2,\ldots,f_{m+1}):\, f_j - f_{j+1} = h_j}
\quad\quad  \prod_{e:\, e \in E}  \gamma(f_1(e),f_2(e),\ldots,f_{m+1}(e)).  
\end{align*}
Summing over all $m$-tuples of flows, we obtain the identity in Lemma 2.3. 
\qed

\vskip 0.4in 
Lemma 2.3 will be used in the following way.  We specialize the indeterminates by choosing a {\sl weight function} 
$\gamma$ defined from $(m+1)$-tuples in $\mathbb{B}^{m+1}$ to a commutative ring $\mathbb{S}.$  The weight function $\gamma$ is extended to $(m+1)$-tuples of functions in the  following way: if $\vec{f} =(f_1,f_2,\ldots,f_{m+1}),$ then  
$$
\gamma(\vec{f}) = \prod_{e:\, e \in E}  \gamma(f_1(e),f_2(e),\ldots,f_{m+1}(e)).
$$    
Let 
$$
\mathcal{F} = \{(f_1,f_2,\ldots,f_{m+1}): \,f_j: E \to \mathbb{B}, \,f_j - f_{j+1}  \,\,\mathrm{are}\,\,\mathrm{flows} \} 
$$
and $\epsilon$ be a element in the ring $\mathbb{S}.$  The {\sl parcel} $\mathcal{P}(\epsilon)$ is the subset of 
$\mathcal{F}$ defined by  
$$
\mathcal{P}(\epsilon) = \{ \vec{f}: \vec{f} \in \mathcal{F}, \, \gamma(\vec{f}) = \epsilon\}.
$$
Parcels are disjoint and their union is $\mathcal{F}.$  Hence,  
$$
\sum_{\epsilon:\, \epsilon \in \mathbb{S}} \epsilon| \mathcal{P}(\epsilon)|
= \sum_{\vec{f}: \, \vec{f} \in \mathcal{F}}  \gamma (\vec{f}).  
\eqno(2)$$ 

The right-hand sum in Eq. (2) is a specialization of the left-hand sum in Lemma 2.3.  Combining Eq. (2) and Lemma 2.3, we have the following lemma.  

\vskip 0.3in\noindent 
{\bf Lemma 2.4.}  Let $G$ be a totally-$\mathbb{P}$ matrix, $\aA$ be a $\mathbb{P}$-module,  $ \mathbb{B} \subseteq \aA,$ and $\gamma: \mathbb{B}^{n+1} \to \mathbb{S}$ be a weight function.  Then
$$
\sum_{\epsilon:\, \epsilon \in \mathbb{S}} \epsilon| \mathcal{P}(\epsilon)|
= 
\sum_{(h_1,h_2,\ldots,h_m): \, h_j \,\text{are flows}}
\quad
\prod_{e:\,e \in E} \,\,\left(
\sum_{(b_1,b_2,\ldots,b_{m+1}): \, b_j \in \mathbb{B},\, b_j - b_{j+1} = h_j(e)}
\gamma (b_1,b_2,\ldots,b_{m+1})
\right).
$$

\vskip 0.2in   
When the weight function is chosen suitably, the left-hand sum in Lemma 2.4 can be written as (a multiple of) an evaluation of the rank generating polynomial.   Whimsically, we might think of this as extracting the grin from the Cheshire cat.   
In this paper, we usually begin with a congruence condition, define parcels with a weight function simulating the congruence condition, 
and then use Lemmas 2.4 and 2.2 to obtain an evaluation of the rank generating polynomial.  The point of evaluation is derived at the end.  

We can also go in reverse in certain cases, starting from a point given beforehand.  We illustrate this by an example.      Suppose 
that we wish to obtain an interpretation for the characteristic polynomial $\chi(M(G);q)$ evaluated at an integer $q$ greater than $1$ using parcels for a totally unimodular rank-$r$ matrix $G.$  By Eq.~(1), this is equivalent to evaluating the rank generating polynomial at $(-q,-1)$ or $(-1,-q).$        
We start with a matrix $H$ (having rank $|E|-r$) orthogonally dual to $G$ and choose the Abelian group to be the integers $\mathbb{Z}_q$ modulo $q$ under addition.  If $\gamma$ is a weight function such that 
$$
\sum_{b,c:\,b-c =0} \gamma (b,c) = q-1  \,\,\, \text{and for}\,\,a \neq 0, \sum_{b,c:\,b-c =a} \gamma (b,c) = -1,  
\eqno(3)$$ 
then by Lemma 2.2, the right-hand side of the equation in Lemma 2.4 (applied to $H$) simplifies in the following way:  
\begin{eqnarray*}
 \sum_{h: \,h \,\text{is  a  flow on}\, H}
(q-1)^{|h^{-1}(0)|}  (-1)^{|E| - |h^{-1}(0)|}   
&=& (-1)^{r} q^{|E|-r} R( M(H);-1,-q) 
\\
&=& 
(-1)^{r} q^{|E| - r} R((M(H)^{\perp};-q,-1) 
\\
&=& q^{|E|-r} \chi(M(G);q). 
\end{eqnarray*} 
Thus, to obtain an interpretation for $\chi(M(G);q),$ it suffices to construct a suitable weight function with a combinatorial or algebraic interpretation.  

Observe that for each $a$ in $\mathbb{Z}_q,$ there are $q$ pairs $(b,c)$ such that 
$b-c =a.$  Thus, when $q$ is odd, we can construct a weight function $\gamma$ satisfying Eq.~(3) in the following way. Begin by setting  
$\gamma(0,0) = 0,$ and for $a \neq 0,$ $\gamma(a,a) = 1.$   For each $a \neq 0,$ choose an even number $s,$  $0 \leq s \leq q-3,$ and $s$ pairs $(b,c)$ such that $b-c=a$ and set $\gamma (b,c) =0.$  Next, choose 
$(q-s-1)/2$ pairs $(b^\prime,c^\prime)$ such that $b^\prime-c^\prime =a$   
and set $\gamma(b,c) = 1.$  Finally, for the remaining $(q-s+1)/2$ pairs $(b^{\prime\prime},c^{\prime\prime})$ such that $b^{\prime\prime} - c^{\prime\prime} =a,$ set $\gamma (b^{\prime\prime},c^{\prime\prime}) = -1.$  
When $q$ is even, a similar construction works, except that we choose an odd number $s.$    

To obtain an explicit example, suppose that $q$ is odd.  If $d$ is an integer modulo $q,$ let $\mathrm{Rem}(d,q)$ be 
the non-negative remainder when $d$ is divided by $q.$  It is easy to check that if $a \not\equiv 0 \,\,\mathrm{mod} \,q,$ then there are $(q+1)/2$ pairs $(b,c)$ such that $b-c =a$ and $\mathrm{Rem}(b+c,q)$ is even.  
Let $\gamma:\, \mathbb{Z}_q \times \mathbb{Z}_q \to \mathbb{C}$ be the weight function defined by $\gamma(0,0) = 0$ and for $(a,b) \neq (0,0),$ 
$$
\gamma(a,b) = \left\{ 
\begin{array}{cc} -1 &\,\text{if } \mathrm{Rem}(a+b,q) \,\,\text{is even}, 
\\
1 &\,\text{if } \mathrm{Rem}(a+b,q) \,\,\text{is odd.} 
\end{array} \right.
$$ 
Then $\gamma$ defines three parcels in the set $\{(f,g):\, f,g:E \to \mathbb{Z}_q, \, f-g \,\,\text{is a flow of}\,\, H\},$  given by   
\begin{align*} 
\mathcal{P}(0) &= \{(f,g): \, \text{for some}\,\, e, \, (f(e),g(e)) = (0,0)\}, 
\\
\mathcal{P}(1) &= \{(f,g): \, \text{for all}\,\, e, \, (f(e),g(e)) \neq (0,0) \,\,\text{and}\,\,  |\{e:\, \mathrm{Rem}(f(e)+g(e),q) \,\text{is even} \}| \equiv 0 \, \mathrm{mod} \, 2 \}, 
\\ 
\mathcal{P}(-1) &= \{(f,g): \, \text{for all}\,\, e, \, (f(e),g(e)) \neq (0,0) \,\,\text{and} \,\,  |\{e:\, \mathrm{Rem}(f(e)+g(e),q) \,\text{is even} \}| \equiv 1 \, \mathrm{mod} \, 2 \}.
\end{align*}
 We can now apply Lemmas 2.4 and 2.2 to obtain the following proposition.         

\vskip 0.2in\noindent 
{\bf Proposition 2.5.}  Let $G$ be a totally unimodular matrix and $q$ be an odd integer greater than $1.$  Then  
$$
|\mathcal{P}(1)| - |\mathcal{P}(-1)|  =  q^{|E|-r}  \chi(M(G);q).  
$$

 \vskip .40in \noindent
{\bf\large    3.  Parcels defined using Hamming distances}

\vskip 0.2in\noindent
 We shall consider two sets of pairs of functions. Let $\aA$ be an Abelian group of order $q$ and $E$ be a set. Let  
\begin{eqnarray*}
\mathcal{F} & = & \{(f,g): \, f,g:E \to \aA, \, f-g \,\,\mathrm{is}\,\,\mathrm{a}\,\,\mathrm{flow}\},  
\\
\mathcal{F^{\times}} & = & \{(f,g):\, f,g: E \to \aA^{\times}, \,f-g \,\,\mathrm{is}\,\,\mathrm{a}\,\,\mathrm{flow}\},   
\end{eqnarray*}  
where $\aA^{\times}$ is the set of nonzero elements in $\aA.$  Since a pair $(f,g)$  in $\mathcal{F}$  defines uniquely a pair $(f,h),$ where $h = f - g$ and $h$ is a flow, and conversely, $|\mathcal{F}|=q^{|E|+r}.$   For the size of $\mathcal{F}^{\times},$  see Corollary 3.4.   

Recall that the {\sl Hamming distance} between two functions $f:E \to \aA$ is the number of elements $e$ in $E$ such that $f(e) \neq g(e),$ or equivalently, the size of $\mathrm{supp} (f-g).$ 

\vskip 0.2in\noindent 
{\bf 3.1.  Parcels in $\mathcal{F}.$}   
\vskip 0.2in  

Let $\sigma$ be an integer greater than $1.$  For $k = 0 ,1, 2, \ldots,\sigma - 1,$ let
\begin{eqnarray*}
\mathcal{O}(k,\sigma) = \{(f,g):\, (f,g) \in \mathcal{F},\,
|\mathrm{supp}(f-g)| \equiv k\,\,\mathrm{mod}\, \sigma\}.  
\end{eqnarray*}
\vskip 0.3in\noindent
{\bf  Theorem 3.1.}
Let $\sigma$ be an integer greater than $1,$  
$\omega$ be a $\sigma$-th primitive root of unity, $G$ be a rank-$r$ totally-$\mathbb{P}$ matrix with column set $E,$  and $\aA$ be a $\mathbb{P}$-module of order $q.$    Then 
\begin{eqnarray*}
\sum_{k=0}^{\sigma - 1}  \omega^k |\mathcal{O}(k,\sigma)|
=
\omega^{|E|-r} (1 - \omega)^r q^{|E|}  R \left( M(G); \frac {q \omega}{1 -\omega} , \frac {1-\omega}{\omega} \right).
\end{eqnarray*}

\vskip 0.1in\noindent 
{\bf Proof.}
 Consider the weight function
$\gamma: \aA \times \aA \to \mathbb{C}$
defined by $\gamma(a,a) = 1$ and $\gamma (a,b) = \omega$ if $a \neq b.$
Then
$$
\sum_{b,c: \, b-c = 0}
\gamma (b,c)
= q, 
\,\,\mathrm{and}\,\, \mathrm{for} \,\, a \neq 0, \,\, 
\sum_{b,c: \, b-c = a}   \gamma (b,c)  = \omega q .
$$
In addition, if $f,g:\, E \to \aA$ are functions, then 
$$
\prod_{e:\,e \in E} \gamma (f(e),g(e)) \,=\,  \prod_{e:\,e \in E, \, f(e) \neq g(e)} \omega 
\,\,=\,\,  \omega^{|\mathrm{supp}(f-g)|}.
$$ 
Hence, by Lemmas 2.4 and 2.2, we have  
\begin{eqnarray*}
\sum_{k=0}^{\sigma - 1}  \omega^k |\mathcal{O}(k,\sigma)|
&=&
\sum_{f,g:\,f - g \,\, \text{is a flow}} \,\,\,\,
\prod_{e:\,e \in E} \gamma (f(e),g(e))
\\
&=& 
\sum_{h:\, h \,\, \text{is a flow}} \,\,  q^{|h^{-1}(0)|} (\omega q)^{|E| - |h^{-1}(0)|}
\\
&=&
\omega^{|E|} q^{|E|} 
 \sum_{h:\, h \, \text{is a flow}} \,\,  \omega^{- |h^{-1}(0)|}
\\
&=&
\omega^{|E|-r} (1 - \omega)^r q^{|E|} R \left( M(G); \frac {q \omega}{1 -\omega} , \frac {1-\omega}{\omega} \right).
\end{eqnarray*}
\qed

\vskip 0.2in\noindent
Theorem 3.1 is the simplest result about parcels.  It can almost be derived by direct substitution from Lemma 2.1.  We shall discuss this further in Section 7.  For now, we note the special case when $\sigma = 2$ (and hence, $\omega = -1$).   

\vskip 0.1in\noindent 
{\bf Corollary 3.2.}  $\displaystyle{\,\, 
|\mathcal{O}(0,2)| - |\mathcal{O}(1,2)| = (-1)^{|E|-r} 2^r q^{|E|}  R \left( M(G); -\frac {q}{2} , -2 \right). } $

\vskip 0.4in\noindent 
{\bf 3.2.  Parcels in $\mathcal{F}^{\times}.$}   
\vskip 0.2in  
  
Let $\sigma$ be a positive integer.  For $k = 0 ,1, 2, \ldots,\sigma - 1,$ let
$$
\mathcal{O}^{\times}(k,\sigma) = \{(f,g):\, (f,g) \in \mathcal{F}^{\times},\,
|\mathrm{supp}(f-g)| \equiv k\,\,\mathrm{mod}\, \sigma\}.  
$$

\vskip 0.2in\noindent
{\bf Theorem 3.3.}
Let $\sigma$ be a positive integer,  $\omega$ be a primitive $\sigma$-th root of unity,  $G$ be a rank-$r$ totally-$\mathbb{P}$ matrix with column set $E,$  and 
$\aA$ be a $\mathbb{P}$-module of order $q,$  where $q \neq 2.$    Then 
$$
\sum_{k=0}^{\sigma - 1}  \omega^k |\mathcal{O}^{\times}(k,\sigma)|
=
\omega^{|E|-r}(q - 2)^{|E|-r}[(1 - \omega)q + 2\omega -1]^r
R \left( M(G); \frac {\omega q(q - 2)}{(1 - \omega)q + 2\omega - 1} , \frac {(1 - \omega)q + 2\omega - 1}{\omega (q - 2)} \right).
$$

\vskip 0.1in\noindent 
{\bf Proof.}
We use the weight function
$\gamma: \aA^{\times} \times \aA^{\times} \to \mathbb{C}$
defined by $\gamma(a,a) = 1$ and $\gamma (a,b) = \omega$ if $a \neq b.$  Then
$$
\sum_{b,c: \, b-c = 0,\,b\neq0,\,c\neq0}
\gamma (b,c)
= q-1, 
$$
and for $a \neq 0,$
$$
\sum_{b,c: \, b-c = a,\,b\neq 0,\,b\neq a}
\gamma (b,c)
= \sum_{b: b \neq 0, \, b \neq a} \gamma (b,b-a)
= \omega (q - 2).
$$
By the argument in the proof of Theorem 3.1, we can finish with the following calculation: 
\begin{eqnarray*}
\sum_{k=0}^{\sigma - 1}  \omega^k |\mathcal{O}^{\times}(k,\sigma)|
&=& 
\sum_{h:\, h \,\, \mathrm{is} \,\,\mathrm{a}\,\,\mathrm{flow}} \,\,  (q-1)^{|h^{-1}(0)|} [\omega (2 - q)]^{|E| - |h^{-1}(0)|}
\\
&=&
\omega^{|E|} (q - 2)^{|E|}
 \sum_{h:\, h \,\, \mathrm{is} \,\,\mathrm{a}\,\,\mathrm{flow}} \,\,  \left( \frac {q-1}{\omega (q - 2)} \right)^{|h^{-1}(0)|}.
\end{eqnarray*}
\qed

\vskip 0.2in \noindent 
The condition that $q \neq 2$ in Theorem 3.3 is needed because we need to divide by $q-2$ in the proof.  
For completeness, we note that when 
$|\aA| = 2,$   there is only one function $E \to \aA^{\times},$ and hence,  $|\mathcal{F}^{\times}| = 1,$  
$|\mathcal{O}^{\times}(0,\sigma)| = 1,$ and for $k \neq 0,$ $|\mathcal{O}^{\times}(k,\sigma)|  = 0.$   

We note three special cases of Theorem 3.3 as corollaries.   The first case is when $\sigma =1$ (so that there is one parcel, $\mathcal{F}^{\times}$ itself).    

\vskip 0.2in\noindent 
{\bf Corollary 3.4.}
If $q \neq 2,$ then 
$$
|\mathcal{F}^{\times}| = (q - 2)^{|E|-r} R \left( M(G); q(q - 2) , \frac {1}{q - 2} \right).
$$

\vskip 0.2in
The second case is when $\sigma =2$ (and $\omega = -1$).  

\vskip 0.2in\noindent 
{\bf Corollary 3.5.}  Let $G$ be a rank-$r$ totally-$\mathbb{P}$ matrix and $\aA$ be a $\mathbb{P}$-module of order $q,$ where $q \neq 2.$   Then  
$$
|\mathcal{O}^{\times}(0,2)| - |\mathcal{O}^{\times}(1,2)|
=
(-1)^{|E|-r}(q - 2)^{|E|-r}(2q - 3)^r
R \left( M(G); - \frac {q(q - 2)}{2q - 3} , - \frac {2q - 3}{q - 2} \right).
$$
In particular, when $G$ is a rank-$r$ ternary matrix and $\aA = \mathrm{GF}(3),$ then 
$$
|\mathcal{O}^{\times}(0,2)| - |\mathcal{O}^{\times}(1,2)| = (-1)^{|E|-r} 3^r R(M(G);-1,-3) 
=
3^r \chi (M(G)^{\perp}; 3).
$$

\vskip 0.2in\noindent 
A concrete example of the second part of Corollary 3.5  is when $\aA = \mathrm{GF}(3)$ and $G$ is the vertex-edge matrix of 
a graph $\Gamma (V,E)$ with $c$ connected components.  If we fix a total order on $V,$ then 
a function from  $E$ to $\{1,-1\}$ defines an orientation on the edges and conversely.  Thus  
$\mathcal{O}^{\times}(k,2)$ is the set of pairs $(f,g)$ of orientations such that $f - g$ is a $\mathrm{GF}(3)$-tension and the number of edges where $f$ and $g$ disagree is congruent to $k$ modulo $2.$  Corollary 3.5 says that    
\begin{eqnarray*}  
|\mathcal{O}^{\times}(0,2)| - |\mathcal{O}^{\times}(1,2)|
= 3^{|V|-c} F(\Gamma;3),
\end{eqnarray*}
where $F(\Gamma;\lambda)$ is the flow polynomial.  This situation can be dualized.  Let $G$ be the cycle-edge matrix of $\Gamma,$  $\bar{\mathcal{O}}^{\times}(k,2)$ be the set of pairs of the orientations such that $f-g$ is a $\mathrm{GF}(3)$-flow, and the number of edges where $f$ and $g$ disagree is congruent to $k$ modulo $2.$ 
Then 
$$
|\bar{\mathcal{O}}^{\times}(0,2)| - |\bar{\mathcal{O}}^{\times}(1,2)|
=
3^{|E|-|V|} P(\Gamma;3),
$$
where $P(\Gamma;\lambda)$ is the chromatic polynomial. 

The third and last case is when $\sigma = 6$ and $q = 3.$         

\vskip 0.2in\noindent 
{\bf Corollary 3.6.}
Let $G$ be a rank-$r$ ternary matrix and $\aA = \mathrm{GF}(3).$  Then 
$$
\sum_{k=0}^5  \ee^{k \pi \iota/3} |\mathcal{O}^{\times}(k,6)| = 
\ee^{|E| \pi \iota/3} (-\sqrt{3} \iota)^r   R(M(G); \sqrt{3} \iota, - \sqrt{3} \iota ).
$$

\vskip 0.1in\noindent 
{\bf Proof.}  Setting $q = 3$ in Theorem 3.3, we obtain 
$$
\sum_{k=0}^5   \omega^k |\mathcal{O}^{\times}(k,6)| = 
\omega^{|E|-r} (2 - \omega)^r  R\left(M(G); \frac {3 \omega}{2 - \omega}, \frac {2-\omega}{\omega} \right).
$$
To finish the proof,  we set $\omega =  \ee^{\pi \iota/3}$ and use the facts $2 - \omega = \sqrt{3} e^{-\pi \iota/6}$ and $\frac {\omega}{2 - \omega} = i/\sqrt{3}.$ 
\qed

\vskip .33in \noindent
{\bf\large  4.  Parcels defined using supports }

\vskip 0.2in\noindent
In this section, we consider parcels in $\mathcal{F}$ defined using supports.   
Let $\sigma$ be an integer greater than $1,$ and $\alpha$ and $\beta$ be integers.
For $k = 0, 1, 2, \ldots, \sigma - 1,$  let
$$
\mathcal{M}_{\alpha,\beta}(k,\sigma) = \{ (f,g): \, (f,g) \in \mathcal{F},  \,\,
\alpha |\mathrm{supp}(f)| + \beta |\mathrm{supp}(g)| \equiv k \,\,\mathrm{mod}\, \sigma\}.
$$

\vskip 0.2in\noindent 
{\bf 4.1.  The general case.}    
\vskip 0.2in 

We begin with the generic theorem. 

\vskip 0.2in\noindent
{\bf  Theorem 4.1.}
Let $\sigma$ be an integer greater then $1,$  $\alpha$ and $\beta$ be integers not congruent to $0$ modulo $\sigma,$ $\omega$ be a primitive $\sigma$-th root of unity,  $G$ be a rank-$r$ totally-$\mathbb{P}$ matrix with column set $E,$ 
and $\aA$ be a $\mathbb{P}$-module of order $q.$   
If  $\omega^{\alpha} + \omega^{\beta} + (q - 2)\omega^{\alpha + \beta} \neq 0,$   then 
\begin{eqnarray*}
\sum_{k=0}^{\sigma - 1} \omega^k |\mathcal{M}_{\alpha,\beta}(k,\sigma)| 
&=& 
(1-\omega^{\alpha})^r(1 - \omega^{\beta})^r [\omega^{\alpha} + \omega^{\beta} + (q - 2)\omega^{\alpha + \beta}]^{|E|-r}
\\
&&  \times \,\,  
R \left( M(G);  \frac {q[\omega^{\alpha} + \omega^{\beta} + (q - 2)\omega^{\alpha + \beta}]}
{(1 - \omega^{\alpha})(1 - \omega^{\beta})}, 
\frac {(1 - \omega^{\alpha})(1 - \omega^{\beta})}
{\omega^{\alpha} + \omega^{\beta} + (q - 2)\omega^{\alpha + \beta}} \right).
\end{eqnarray*}

\vskip 0.1in\noindent
{\bf Proof.}   We use the weight function $\gamma: \aA \times \aA \to \mathbb{C}$ 
defined by $\gamma(0,0) = 1,$ and for $a$ and $b$ nonzero elements (not necessarily distinct) of $\aA,$
$$\gamma(a,0) = \omega^{\alpha}, \,\, \gamma(0,a) = \omega^{\beta}, \,\,  \gamma(a,b) = \omega^{\alpha + \beta}.
$$  
Let $f,g:\, E \to \aA$ be functions,  $A = \mathrm{supp}(f),$ and $B = \mathrm{supp}(g). $
Then
$$
\prod_{e:\,e \in E} \gamma (f(e),g(e)) = \omega^{\alpha (|A| - |A \cap B|)} \omega^{\beta (|B| - |A \cap B|)}
\omega^{(\alpha + \beta)|A \cap B|}
=\omega^{\alpha |A| + \beta |B|}.
$$  
Next, observe that 
$$
\sum_{b,c: \, b-c = 0}
\gamma (b,c)  = 1 + (q-1)\omega^{\alpha + \beta},
$$
and for $a \neq 0,$
\begin{eqnarray*}
\sum_{b,c: \, b-c = a}
\gamma (b,c)
&=&  
\gamma (a,0) + \gamma (0,-a) + \sum_{b: b \neq 0, b \neq a} \gamma (b,b-a)
\\
&=&  \omega^{\alpha} + \omega^{\beta} + (q - 2)\omega^{\alpha + \beta}. 
\end{eqnarray*}
Then, as in the proof of Theorem 3.1, we can finish with the following calculations:    
\begin{eqnarray*}
&&  \sum_{k=0}^{\sigma - 1} \omega^k |\mathcal{M}_{\alpha,\beta}(k,\sigma)| 
\\
&=& 
\sum_{h:\,  h \,\mathrm{is} \,\,\mathrm{a}\,\, \mathrm{flow}}  [1 + (q-1)\omega^{\alpha + \beta}]^{|h^{-1}(0)|}
[\omega^{\alpha} + \omega^{\beta} + (q - 2)\omega^{\alpha + \beta}]^{|E| - |h^{-1}(0)|}
\\
&=&
[\omega^{\alpha} + \omega^{\beta} + (q - 2)\omega^{\alpha + \beta}]^{|E|}
\sum_{h:\, h \,\mathrm{is} \,\,\mathrm{a}\,\, \mathrm{flow}}
\left( \frac {1 + (q-1)\omega^{\alpha + \beta}}{\omega^{\alpha} + \omega^{\beta} + (q - 2)\omega^{\alpha + \beta}} \right)^{|h^{-1}(0)|},
\end{eqnarray*}
and 
$$
 \frac {1 + (q-1)\omega^{\alpha + \beta}}{\omega^{\alpha} + \omega^{\beta} + (q - 2)\omega^{\alpha + \beta}} - 1 = 
\frac {(1 - \omega^{\alpha})(1 - \omega^{\beta})}
{\omega^{\alpha} + \omega^{\beta} + (q - 2)\omega^{\alpha + \beta}}.
$$
\qed

\vskip 0.3in\noindent 
{\bf 4.2.  The special case when $\alpha = \beta = 1.$ }    
\vskip 0.2in 

Let $\alpha = \beta = 1.$  Then 
$$
\mathcal{M}_{1,1}(k,\sigma)  = \{(f,g):\, (f, g) \in \mathcal{F},\,\, 
|\mathrm{supp}(f)| + |\mathrm{supp}(g)| \equiv k \,\,\mathrm{mod}\,\sigma \}
$$
and we have the following specialization of Theorem 4.1.

\vskip 0.2in\noindent
{\bf Theorem 4.2.}
Let $\sigma$ be an integer greater than $1,$ $\omega$ be a primitive $\sigma$-th root of unity, $G$ be a rank-$r$ totally-$\mathbb{P}$ matrix with column set $E,$ and $\aA$ be a $\mathbb{P}$-module of order $q.$  Then except when $\sigma = 2$ and $q=4,$  
$$
\sum_{k=0}^{\sigma - 1} \omega^k |\mathcal{M}_{1,1}(k,\sigma)|
=
(1 - \omega)^{2r} (2\omega + \omega^2(q-2))^{|E|-r}
R \left( M(G);  \frac {q (2 \omega  + \omega^2(q-2))}{(1 - \omega)^2}, \frac {(1 - \omega)^2}{2\omega + \omega^2(q-2)} \right).
$$

\vskip 0.2in \noindent  
The exception in Theorem 4.2 is when $\sigma = 2$ and $q=4.$  When this happens, $2\omega + \omega^2 (q-2) = 0.$ 
Following the proof of Theorem 4.1 with $\alpha = \beta = 1$ up to the point when we need to divide by $2 \omega + \omega^2(q-2),$  we obtain   
$$
|\mathcal{M}_{1,1}(0,2)| - |\mathcal{M}_{1,1}(1,2)| = 
\sum_{h:\, h \,\mathrm{is}\,\mathrm{a} \,\mathrm{flow}} 4^{|h^{-1}(0)|}0^{|E|- |h^{-1}(0)|} = 4^{|E|}.  
$$

There are several corollaries of Theorem 4.2 giving evaluations of the rank generating polynomial at rational or real numbers.  

\vskip 0.2in\noindent 
{\bf Corollary 4.3.} 
Let $\sigma = 2$ and $q \neq 4.$  Then  
$$
|\mathcal{M}_{1,1}(0,2)| - |\mathcal{M}_{1,1}(1,2)| = 
4^r (q -4)^{|E|-r} R \left( M(G);  \frac {q(q-4)}{4}, \frac {4}{q-4} \right).
$$
In particular, when $G$ is a binary matrix and $\aA = \mathrm{GF}(2),$    
$$
|\mathcal{M}_{1,1}(0,2)| - |\mathcal{M}_{1,1}(1,2)| = 
2^{|E|+r} (-1)^{|E|-r}R ( M(G);  -1, -2 ) = 2^{|E|+r} \chi (M(G)^{\perp};2).
\eqno(4)$$

\vskip 0.2in\noindent 
Eq.~(4) in Corollary 4.3 is an analog (extended from graphs to binary matroids) 
of Proposition 12 in \cite{Goodall}.  
It may be worthwhile to give a combinatorial proof of Eq.~(4).  This will give an independent verification of the evaluation and tell us what counting information that evaluation contains.  We need two easy facts in the proof.  First, if $h = f -g,$  then  
 $$
|\mathrm{supp}(h)|= |\mathrm{supp}(f)| +|\mathrm{supp}(g)| - 2|\mathrm{supp}(f) \cap \mathrm{supp}(g)| 
\equiv |\mathrm{supp}(f)| +|\mathrm{supp}(g)|\,\, \mathrm{mod} \,\,2.
$$
Second, for each flow $h,$ there are $2^{|E|}$ pairs $(f,g)$ such that $f-g = h.$   These two facts imply that 
$|\mathcal{M}_{1,1}(0,2)|$  (respectively, $|\mathcal{M}_{1,1}(1,2)|$) equals $2^{|E|}$ times the number of flows 
with support of even (respectively, odd) size.  
Since flows with support of even size form a subspace, either 
(a)  all flows have support of even size, or (b) half the flows have support of even size, and the other half have support of odd size.  In case (a), $|\mathcal{M}_{1,1}(0,2)| = 2^{|E|} 2^r$ and $|\mathcal{M}_{1,1}(1,2)| = 0.$  Moreover, all cocircuits of the matroid $M(G)$ are even and hence, the dual $M(G)^{\perp}$ is $\mathrm{GF}(2)$-affine.  
By Lemma 1.4, $\chi(M(G)^{\perp};2)=1$ and thus, the two sides of Eq.~(4) are equal.  
In case (b),  
$|\mathcal{M}_{1,1}(0,2)| = |\mathcal{M}_{1,1}(1,2)|,$ and the left-hand side of Eq.~(4) is $0.$  Since supports of flows of binary matrices are disjoint union of cocircuits, $M(G)$ has a cocircuit of odd size.  Therefore,  $M(G)^{\perp}$ is not $\mathrm{GF}(2)$-affine, $\chi (M(G)^{\perp};2) = 0,$ and both sides of Eq.~(4) equals 0.  This completes the combinatorial proof of Eq.~(4).   

When $q=2,$ Theorem 4.2 gives evaluations of the rank generating polynomial at real values.  

\vskip 0.2in\noindent 
{\bf Corollary 4.4.} 
Let $G$ be a binary matrix, $\aA = \mathrm{GF}(2),$  $\theta = 2 \rho \pi /\sigma,$ $\omega = \ee^{\iota \theta},$  and $\rho$ be an integer relatively prime to $\sigma.$   Then 
$$
\sum_{k=0}^{\sigma - 1} \omega^k |\mathcal{M}_{1,1}(k,\sigma)|
=
(\omega - 1)^{2r} (2\omega )^{|E|-r}
R \left( M(G);  \frac {2}{\cos \theta - 1}, \cos \theta - 1 \right). 
$$

\vskip 0.2in\noindent 
The cases $\rho = 1$ and $\sigma = 2,3,4,$ or $ 6$ give interpretations
of the rank generating polynomial at points with rational coordinates. The case $\sigma = 2$ was covered in Eq.~(4).  
The cases $\sigma = 3$ ($\cos \theta =  -1/2$),  $\sigma = 4$   ($\cos \theta =  0$), and $\sigma = 6$  ($\cos \theta =  1/2$). 

\vskip 0.2in\noindent 
{\bf Corollary 4.5.}  Let $G$ be a binary matrix and $\aA = \mathrm{GF}(2).$   Then 
\begin{eqnarray*} 
&& 
\sum_{k=0}^{2} \ee^{2k \pi \iota /3} |\mathcal{M}_{1,1}(k,3)|  
=
\ee^{2 |E| \pi \iota /3} 2^{|E|-r} (-3)^r  R \left( M(G);  - \frac {4}{3}, - \frac {3}{2} \right),  
\\
&&
 \sum_{k=0}^3  \iota^k |\mathcal{M}_{1,1}(k,4)| 
= (-1)^r (2 \iota)^{|E|} R(M(G); -2,-1) = (2\iota)^{|E|} \chi ( M(G);  2), 
\\
&& 
\sum_{k=0}^5   \ee^{k \pi \iota /3} |\mathcal{M}_{1,1}(k,6)|
=
(-1)^r \ee^{|E| \pi \iota /3} 2^{|E|-r}  R \left( M(G);  - 4, - \frac {1}{2} \right).
\end{eqnarray*}

\vskip 0.2in\noindent 
This first equation in Corollary 4.5, when $\sigma = 3,$ is an analog of Theorem 15 in \cite{Goodall}.   We note one more case, when $\sigma = q = 3.$

\vskip 0.2in\noindent 
{\bf Corollary 4.6.}  
 If $G$ is a rank-$r$ ternary matrix and $\aA = \mathrm{GF}(3),$ then 
$$
|\mathcal{M}_{1,1}(0,3)| + \ee^{2 \pi \iota /3} |\mathcal{M}_{1,1}(1,3)| + \ee^{-2 \pi \iota /3} |\mathcal{M}_{1,1}(2,3)|    = 
 (\sqrt{3} \ee^{5\pi \iota/6})^{|E| + r}  R( M(G);  \sqrt{3}  \ee^{ -5\pi \iota/6},  \sqrt{3} \ee^{5 \pi \iota/6}).
$$

\vskip 0.1in\noindent 
{\bf Proof.} When $q=3$ and $\omega =  \ee^{2 \pi \iota/3},$  $2\omega + \omega^2(q-2) = \omega - 1.$   
Theorem 4.3 yields 
$$
|\mathcal{M}_{1,1}(0,3)| + \ee^{2 \pi \iota /3} |\mathcal{M}_{1,1}(1,3)| + \ee^{-2 \pi \iota /3} |\mathcal{M}_{1,1}(2,3)|    = (-1)^{|E| - r} (1-\omega)^{|E|+r} R\left( M(G);  \frac {3}{\omega - 1}, \omega - 1   \right). 
$$
To finish the proof, use $\omega - 1 = \sqrt{3}  \ee^{5 \pi \iota/6}.$   Note that 
$ R ( M(G);  3 /(\omega - 1), \omega - 1) = T(M(G);\omega^2,\omega);$ thus, we have an evaluation of the Tutte polynomial at $(\omega^2,\omega).$
\qed

\vskip 0.4in\noindent 
{\bf 4.3.  The special case when $\alpha = 1$ and $\beta = -1.$ }    
\vskip 0.2in 

Let $\alpha = 1$  and $\beta = -1.$   Then 
$$
\mathcal{M}_{1,-1}(k,\sigma) = \{ (f,g): \,(f,g) \in \mathcal{F},\,\,
|\mathrm{supp}(f)| - |\mathrm{supp}(g)| \equiv k \,\,\mathrm{mod}\, \sigma\}.
$$
\vskip 0.2in\noindent
{\bf Theorem 4.7.}
Let $\sigma$ be an integer greater than $1,$  $\theta = 2\rho \pi/\sigma,$ where $\rho$ is relatively prime to $\sigma,$   $\omega = \ee^{\iota \theta},$  $G$ be a rank-$r$ totally-$\mathbb{P}$ matrix with column set $E,$ 
and $\aA$ be a $\mathbb{P}$-module of order $q.$   With three exceptions, $\sigma = 2$ and $q = 4,$  $\sigma =3$ and $q = 3,$ or $\sigma = 4$ and $q = 2,$ we have 
$$
\sum_{k=0}^{\sigma - 1} \omega^k |\mathcal{M}_{1,-1}(k,\sigma)|
=
(2 - 2 \cos \theta)^r (2 \cos \theta  - 2 + q)^{|E|-r}
R \left( M(G);  \frac {q(2 \cos \theta - 2 + q)}{2 - 2\cos \theta}, \frac {2 - 2\cos \theta}{2 \cos \theta - 2 + q} \right).
$$

\vskip 0.1in\noindent
{\bf Proof.}  Set $\alpha = 1$ and $\beta = - 1$ in Theorem 4.1 and observe that $\omega + \omega^{-1} = 2\cos \theta$  and $\omega\omega^{-1} = 1.$  
\qed 

\vskip 0.2in \noindent
The three exceptional cases in Theorem 4.7 occur because $2 \cos \theta - 2 + q = 0.$   
We note, for completeness, that  
\begin{eqnarray*}  
&& |\mathcal{M}_{1,-1}(0,2)| - |\mathcal{M}_{1,-1}(1,2)| =  4^{|E|},  
\\
&& |\mathcal{M}_{1,-1}(0,3)| + \ee^{2 \pi \iota/3} |\mathcal{M}_{1,-1}(1,3)| + \ee^{-2 \pi \iota/3} |\mathcal{M}_{1,-1}(2,3)| 
= 
3^{|E|},  
\\
&& |\mathcal{M}_{1,-1}(0,4)| + \iota |\mathcal{M}_{1,-1}(1,4)| - |\mathcal{M}_{1,-1}(2,4)| - \iota |\mathcal{M}_{1,-1}(3,4)|  =
 2^{|E|}. 
\end{eqnarray*}

Theorem 4.7 says that a sum of complex numbers equals a real number. Thus, it gives
two equations between real numbers.  

\vskip 0.2in \noindent 
{\bf Corollary 4.8.}           
\begin{eqnarray*}
&& \sum_{k=0}^{\sigma - 1} (\cos k \theta) |\mathcal{M}_{1,-1}(k,\sigma)|
\\
&=&
(2-2 \cos \theta)^r (2 \cos \theta - 2 + q)^{|E|-r}
R \left( M(G);  \frac {q (2 \cos \theta - 2 + q)}{2 - 2 \cos \theta},
\frac {2 - 2\cos \theta}{2 \cos \theta - 2 + q} \right)
\end{eqnarray*} 
and 
$$
\sum_{k=1}^{\sigma - 1} (\sin k \theta) |\mathcal{M}_{1,-1}(k,\sigma)|=0.
\qquad\qquad\qquad \qquad\qquad\qquad\qquad\qquad\qquad$$

\vskip 0.4in\noindent 
The second equation is not unexpected.  Since $|B| - |C| = -[|C| - |B|],$ the {\sl twist} sending  
$(f,g)$ to $(g,f)$ is a bijection from $\mathcal{M}_{1,-1}(k,\sigma)$ to $\mathcal{M}_{1,-1}(\sigma - k,\sigma).$
Hence, we have the following strengthening of the second equation in Corollary 4.8.      

\vskip 0.2in \noindent 
{\bf Proposition 4.9.} $ |\mathcal{M}_{1,-1}(k,\sigma)| = |\mathcal{M}_{1,-1}(\sigma-k,\sigma)|. $ 
  
\vskip 0.2in 
As in subsection 4.2,  the cases $\rho = 1$ and $\sigma = 2,3,4,$ or $ 6$
gives evaluation of the rank generating polynomial at rational numbers. The case $\sigma = 2$ and $q \neq 4$ was done  earlier in Eq.~(4) in Corollary 4.3.  

When $\sigma = 3,$  
$\cos \theta = \cos 2\pi/3 = -\frac {1}{2}.$  By Proposition 4.9, $|\mathcal{M}_{1,-1}(1,3)| =  |\mathcal{M}_{1,-1}(2,3)|.$ Hence,  we have the following special case of Corollary 4.8. 

\vskip 0.2in \noindent 
{\bf Corollary 4.10.}   When $q \neq 3,$
\begin{eqnarray*}
|\mathcal{M}_{1,-1}(0,3)| - |\mathcal{M}_{1,-1}(1,3)|
= 3^r (q -3)^{|E|-r} R \left( M(G);  \frac {q (q - 3)}{3}, \frac {3}{q - 3} \right). 
\end{eqnarray*}
In particular, when $G$ is a binary matrix and $\aA = \mathrm{GF}(2),$   
$$
|\mathcal{M}_{1,-1}(0,3)| - |\mathcal{M}_{1,-1}(1,3)| =
(-1)^{|E|-r} 3^r R \left( M(G);  -\frac {2}{3}, -3 \right).    
$$

\vskip 0.2in 
The second equation in Corollary 4.10 is an analog of Theorem 14 in \cite{Goodall}.
Next, we consider the case when $\sigma = 4,$ $\cos \theta = \cos \pi/2 = 0.$   

\vskip 0.2in \noindent 
{\bf Corollary 4.11.}  When $q \neq 2,$    
$$
|\mathcal{M}_{1,-1}(0,4)| - |\mathcal{M}_{1,-1}(2,4)| = 2^r(q-2)^{|E|-r}  R\left(M(G); \frac {q(q-2)}{2}, \frac {2}{q-2}\right).   
$$

\vskip 0.2in \noindent 
Finally, we consider the case when $\sigma = 6$ and $\cos \theta = \cos \pi/3 = \frac {1}{2}.$    

\vskip 0.2in \noindent 
{\bf Corollary 4.12.}     
\begin{eqnarray*}
 |\mathcal{M}_{1,-1}(0,6)| + |\mathcal{M}_{1,-1}(1,6)| - |\mathcal{M}_{1,-1}(2,6)| - |\mathcal{M}_{1,-1}(3,6)| 
= (q-1)^{|E|-r}  R\left(M(G); q(q-1), \frac {1}{q-1}\right).   
\end{eqnarray*}   

\vskip 0.2in\noindent 
{\bf 4.4.  The special case when $\sigma = 2 \tau,$ $\alpha =1,$ and $ \beta = \tau -1.$ }    
\vskip 0.2in 

In this subsection, we consider the case when $\sigma = 2\tau,$ where $\tau$ is a positive integer, $\alpha = 1,$ and $\beta = \tau - 1.$  
In this case, $2 \tau$ parcels in $\mathcal{F}$ are defined.  For $k = 0, 1, 2, \ldots, 2\tau - 1,$ let
$$
\mathcal{M}_{1,\tau-1}(k, 2\tau) = \{ (f,g): \, (f, g) \in \mathcal{F},\,\, 
|\mathrm{supp}(f)| + (\tau - 1)|\mathrm{supp}(g)| \equiv k \,\,\mathrm{mod}\, 2\tau \}.
$$

\vskip 0.2in\noindent
{\bf Theorem 4.13.}  Let $\tau$ be a positive integer, $\theta = \rho \pi / \tau,$ where $\rho$ is relatively prime to $2\tau,$  
$\omega = \ee^{\iota \theta},$   $G$ be a rank-$r$ totally-$\mathbb{P}$ matrix with column set $E,$  and $\aA$ be a $\mathbb{P}$-module of order $q.$  Then 
\begin{eqnarray*}
&& \sum_{k=0}^{2\tau - 1} \omega^k |\mathcal{M}_{1,\tau-1}(k, 2\tau)|
\\
&&=
(-2 \iota \sin \theta)^r (2 \iota \sin \theta + 2 - q)^{|E|-r}
R \left( M(G);  \frac {-q (2 \iota \sin \theta + 2 - q)}{2 \iota \sin \theta},
\frac {-2 \iota \sin \theta}{2 \iota \sin \theta + 2 -q} \right).
\end{eqnarray*}

\vskip 0.1in\noindent
{\bf  Proof.}   Set $\alpha = 1$ and $\beta = \tau-1$ in Theorem 4.2 and observe that $\omega + \omega^{\tau - 1} = \omega - \omega^{- 1} = 2 \iota \sin \theta$ and $\omega \omega^{\tau - 1} = -1.$ 
\qed
 
\vskip 0.2in              
A typical example of Theorem 4.13 is when $\tau = 6,$  $\omega = \ee^{\pi \iota/6},$
and $\sin \theta = 1/2.$  The congruence condition for defining parcels is 
$ |\mathrm{supp}(f)| + 5 |\mathrm{supp}(g)| \equiv k \,\,\,\mathrm{mod}\,12 $  and we have  
\begin{eqnarray*}
\sum_{k=0}^{11} \ee^{k \pi \iota/6} |\mathcal{M}_{1,5}(k,12)|
=
(- \iota )^r ( \iota + 2 - q)^{|E|-r}
R \left( M(G);  q \iota ( \iota + 2 - q), \frac {- \iota }{ \iota + 2 -q} \right).
\end{eqnarray*}

When $G$ is a rank-$r$ binary matrix and $\aA = \mathrm{GF}(2),$ we obtain evaluations of the characteristic polynomials.  For binary matroids, we may replace a function by its support.  Specifically,  we have  
$$
\mathcal{M}_{1,\tau-1}(k, 2 \tau) =
\{(A,B): A \Delta B \,\,\mathrm{is} \,\, \mathrm{a} \,\,\mathrm{union} \,\,\mathrm{of} \,\,\mathrm{circuits}, \, 
|A| + (\tau-1)|B| \equiv k \,\,\mathrm{mod} \,\, 2 \tau  \},
$$
where $\Delta$ is symmetric difference of sets.  Since $q-2 =0,$ cancellations occur in Theorem 4.13, and  we have 
the following result.  

\vskip 0.2in\noindent 
{\bf Corollary 4.14.}  Let $G$ be a binary rank-$r$ matrix with column set $E$ and $\aA = \mathrm{GF}(2).$   Then 
\begin{eqnarray*}
\sum_{k=0}^{2 \tau - 1} \ee^{k \pi \iota/\tau} |\mathcal{M}_{1,\tau-1}(k,2 \tau)|
&=&
(- 1 )^r (2 \iota \sin \theta)^{|E|}
R ( M(G);  -2,-1) 
\\
&=&  (2 \iota \sin \theta)^{|E|} \chi (M(G);2).
\end{eqnarray*}

\vskip 0.2in\noindent 
Narrowing further to the case $\tau = 2$ (so that the congruence condition is $|A|+|B| \equiv k \,\mathrm{mod}\,4$), we have another result. 

\vskip 0.2in\noindent 
{\bf Corollary 4.15.} Let $G$ be a binary rank-$r$ matrix with column set $E.$  Then 
$$
|\mathcal{M}_{1,1}(|E|,4)| - |\mathcal{M}_{1,1}( |E| + 2,4 )|= 2^{|E|} \chi (M(G);2),
$$
and 
$$
|\mathcal{M}_{1,1}(|E|+1,4)| = |\mathcal{M}_{1,1}( |E| + 3, 4)| 
$$
where the integers $|E|,$  $|E|+1,$ $|E|+2,$ and $|E|+3$ are reduced modulo $4.$  
The first equation implies that the binary matroid $M(G)$ is affine if and only if 
$|\mathcal{M}_{1,1}(|E|,4)| \neq |\mathcal{M}_{1,1}( |E| + 2,4)|.$

\vskip 0.1in\noindent
{\bf Proof.}  
Setting $\tau = 2$ in Corollary 4.14, we have 
$$
|\mathcal{M}_{1,1}(0,4)| + \iota |\mathcal{M}_{1,1}(1,4)| - |\mathcal{M}_{1,1}(2,4)| - \iota|\mathcal{M}_{1,1}(3,4)| = 
(2\iota)^{|E|} \chi (M(G);2).  
$$
Shifting the first parameter $k$ of the parcels by $|E|$ and equating real and imaginary parts, we obtain the two equations in the corollary.  
\qed 

\vskip 0.2in \noindent 
We end with two remarks.  First, Corollary 4.15 can also be obtained as the special case $\sigma = 4$ and $q=2$ of Theorem 4.2.   Second, the case when $G$ is the vertex-edge or cycle-edge matrix of a graph of Corollary 4.15, together with Lemma 1.5,  yield Theorem 1.1. 

\vskip 0.2in\noindent 
{\bf 4.5.  Parcels defined using a binary set-operation on supports. }    
\vskip 0.2in 
 
Let $\circ$ be a binary operation on sets and $\sigma$ be a positive integer.  For $k = 0, 1, 2, \ldots, \sigma - 1,$  let
\begin{eqnarray*} 
\mathcal{M}^{\circ}(k,\sigma)  
=  
\{ (f,g): \, (f, g) \in \mathcal{F}, \,\, |\mathrm{supp}(f) \circ \mathrm{supp}(g)| \equiv k \,\,\mathrm{mod}\, \sigma\}
\end{eqnarray*}
There are sixteen binary operations on subsets of a set $E.$  Familiar examples are union, intersection, and symmetric difference.  Exotic examples are the {\sl Sheffer stroke} $|,$ defined by 
$A\,|\,B = (E \backslash A) \cap (E \backslash B) = E \backslash (A \cup B),$ and 
{\sl implication} $\to,$ defined by $A \to B = (E \backslash A) \cup B.$   

\vskip 0.2in\noindent
{\bf  Theorem 4.16.}
Let $\sigma$ be a integer greater than $1,$ $\omega$ be a primitive $\sigma$-th root of unity,  $G$ be a rank-$r$ totally-$\mathbb{P}$ matrix with column set $E,$ and $\aA$ be a $\mathbb{P}$-module of order $q.$   

\vskip 0.2in \noindent
(a) $\displaystyle{ 
\sum_{k=0}^{\sigma - 1} \omega^k |\mathcal{M}^{\cup}(k,\sigma)| 
= 
(1-\omega)^r (\omega q)^{|E|-r}  
R \left( M(G);  \frac {\omega q^2}{1 -\omega}, \frac {1 - \omega}{\omega q} \right).}
$

\vskip 0.2in \noindent
(b)   If  $2 + \omega(q-2) \neq 0,$ then
$$ 
\sum_{k=0}^{\sigma - 1} \omega^k |\mathcal{M}^{\cap}(k,\sigma)| 
= 
(\omega-1)^r (2 + \omega (q-2))^{|E|-r}  
R \left( M(G);  \frac {q(2 + \omega (q-2))}{\omega - 1},  \frac {\omega - 1}{2 + \omega(q-2)} \right).
$$

\vskip 0.2in \noindent
(c) 
If $2 \omega + q - 2 \neq 0,$ then  
$$
\sum_{k=0}^{\sigma - 1} \omega^k |\mathcal{M}^{\Delta}(k,\sigma)| 
= 
(2- 2 \omega)^r (2 \omega + q - 2)^{|E|-r}  
R \left( M(G);  \frac {q(2 \omega + q - 2)}{2 -2 \omega}, \frac {2 - 2 \omega}{2 \omega + q - 2} \right). 
$$

\vskip 0.2in \noindent
(d) 
$\displaystyle{
\sum_{k=0}^{\sigma - 1} \omega^k |\mathcal{M}^{|}(k,\sigma)| 
= 
(\omega - 1)^r q^{|E|-r}  
R \left( M(G);  \frac {q^2}{\omega -1}, \frac {\omega -1}{q} \right). 
}$

\vskip 0.2in \noindent
(e) 
If $1  + \omega (q - 1) \neq 0,$ then  
$$
\sum_{k=0}^{\sigma - 1} \omega^k |\mathcal{M}^{\to}(k,\sigma)| 
= 
(\omega - 1)^r (1  + \omega (q - 1))^{|E|-r}  
R \left( M(G);  \frac {q(1 + \omega (q - 1))}{\omega -1}, 
\frac {\omega - 1}{1 +  \omega (q - 1)} \right). 
$$

\vskip 0.1in\noindent
{\bf Proof.}    
As examples, we prove (c) and (d).  To prove (c), we use the weight function defined by  
$\gamma (0,0)=\gamma (a,a)=1,$ $\gamma(a,0)=\gamma(0,a) = \omega,$ $\gamma (a,b)=1.$ For this choice of $\gamma,$ 
$$
x - 1 = \frac {q}{2\omega + q - 2} - 1 = \frac {2 - 2 \omega}{2 \omega + q - 2} .
$$
To prove (d), we use the weight function defined by $\gamma (0,0)= \omega,$  $\gamma (a,a)= 1,$ $\gamma(a,0)=\gamma(0,a) = \gamma (a,b)= 1.$   For this choice of $\gamma,$ 
$  x - 1 =  [(\omega + q - 1)/q] - 1 = (\omega-1)/q. $
\qed

\vskip 0.3in 
When $\sigma = 2$ and $\omega = -1,$ we have the following special cases.  

\vskip 0.2in\noindent 
 {\bf  Corollary 4.17.}   Let $G$ be a rank-$r$ totally-$\mathbb{P}$ matrix with column set $E,$ and $\aA$ be a $\mathbb{P}$-module of order $q.$  Then,      
   
\begin{eqnarray*}
\\
|\mathcal{M}^{\cup}(0,2)| - |\mathcal{M}^{\cup}(1,2)| 
&=& 
(-1)^{|E|} (-2)^r  q^{|E|-r}  
R \left( M(G);  - \frac {q^2}{2}, - \frac {2}{q} \right), 
\\
|\mathcal{M}^{|}(0,2)| - |\mathcal{M}^{|}(1,2)| 
&=&  
(-2)^r q^{|E|-r}  
R \left( M(G);  - \frac {q^2}{2}, -\frac {2}{q} \right). 
\end{eqnarray*}  
If $q \neq 4,$ then 
\begin{eqnarray*}
|\mathcal{M}^{\cap}(0,2)| - |\mathcal{M}^{\cap}(1,2)| 
&=& 
(-2)^r (4 - q)^{|E|-r}  
R \left( M(G);  -\frac {q(4-q)}{2},  -\frac {2}{4-q} \right),
\\
|\mathcal{M}^{\Delta}(0,2)| - |\mathcal{M}^{\Delta}(1,2)| 
&=& 
4^r (q - 4)^{|E|-r}  
R \left( M(G);  \frac {q(q - 4)}{4}, \frac {4}{q - 4} \right). 
\end{eqnarray*} 
If $q \neq 2,$ then 
$$
|\mathcal{M}^{\to}(0,2)| - |\mathcal{M}^{\to}(1,2)| 
=   
(-2)^r (2 - q)^{|E|-r}  
R \left( M(G);   \frac {q(q - 2)}{2},  \frac {2}{q - 2} \right). 
$$

\vskip 0.2in 
We note a last special case.  

\vskip 0.2in\noindent 
 {\bf  Corollary 4.18.} 
Let $G$ is a ternary matroid and $\aA = \mathrm{GF}(3).$  Then 
$$
|\mathcal{M}^{\to}(0,3)| + \ee^{2 \pi \iota/3} |\mathcal{M}^{\to}(1,3)| + \ee^{4 \pi \iota/3} |\mathcal{M}^{\to}(2,3)| 
=  
\sqrt{3}^{|E|} \ee^{(3|E| + 5r) \pi \iota /6}   
R \left( M(G);   3 \ee^{- \pi \iota/3},  \ee^{\pi \iota/3} \right). 
$$

\vskip .33in \noindent
{\bf\large   5.  Multivariate versions }

\vskip .20in \noindent
So far, we have been working with pairs of functions.  In this section, we explore parcels of $m$-tuples of functions.  
Let $\aA$ be an Abelian group of order $q.$   If $(a_1,a_2, \ldots,a_m)$ is an $m$-tuple in $\aA^m,$  then there are 
$|\aA|$ $(m+1)$-tuples $(b_1, b_2, \ldots ,b_m, b_{m+1})$ such that for $j = 1,2, \ldots,m,$ $b_j - b_{j+1} = a_j.$  We say that these $(m+1)$-tuples are {\sl associated} with $(a_1,a_2,\ldots,a_m).$   Explicitly, the $(m+1)$-tuples 
$$
(a_1+a_2+\cdots+a_{m-1} +a_{m}+b, a_2+\cdots+a_{m-1} +a_{m}+b, \ldots ,a_{m-1}+a_m + b, a_m + b, b),
$$
where $b$ ranges over $\aA,$ are associated with $(a_1,a_2,\ldots,a_m).$  

 We begin with parcels of triples of functions defined by sums of supports of functions.   
Let $\mathcal{T}(k,\sigma)$ be the parcels 
$$
\{(f_1,f_2,f_3): \, f_j:E \to \aA, \, f_j-f_{j+1} \,\text{are flows},\,  
|\mathrm{supp}(f_1)| +|\mathrm{supp}(f_2)| +|\mathrm{supp}(f_3)| \equiv k \,\text{mod}\,\sigma\}.
$$

\vskip 0.2in \noindent 
{\bf Theorem 5.1.}  Let $\sigma$ be a positive integer, $\omega$ be a primitive $\sigma$-th root of unity, 
$\mathbb{P}$ be a partial field in which $2 \neq 0,$ $G$ be a
rank-$r$ totally-$\mathbb{P}$ matrix, and $\aA$ be a $\mathbb{P}$-module of order $q.$  
Then 
\begin{eqnarray*}
&& \sum_{k=0}^{\sigma-1} \omega^k |\mathcal{T} (k,\sigma)|
\\
&=&
\sum_{
\genfrac{}{}{0pt}{}
{\vec{h}: \, \vec{h} =(h_1,h_2),}{h_1, \, h_2 \,\mathrm{are}\,\,\mathrm{flows}}} 
(1 + (q-1)\omega^3)^{|\tilde{h}^{-1}(0,0)|}
\prod_{a:\, a \neq 0} (\omega + \omega^2 + (q-2) \omega^3)^{|\vec{h}^{-1}(a,0)|+ |\vec{h}^{-1}(0,a)|}  
\\
&& \qquad\qquad\qquad\qquad\qquad 
\times \prod_{a,b:\, a \neq 0, \,b \neq 0} (3\omega^2 + (q-3) \omega^3)^{|\vec{h}^{-1}(a,b)|}.  
\end{eqnarray*}

\vskip 0.1in\noindent 
{\bf Proof.}   We use the weight function defined by  
$\gamma(0,0,0) = 1$ and for $a,b,c$ nonzero elements (not necessarily distinct) of $\aA,$ 
$$
\gamma(a,0,0)=\gamma(0,a,0) =\gamma(0,0,a) =\omega, \,\,\gamma(a,b,0)=\gamma(a,0,b) =\gamma(0,a,b) =\omega^2, \,\, \gamma(a,b,c)=\omega^3. 
$$
Then, for $f_1,f_2,f_3:\,E \to \aA$ with $\mathrm{supp}(f_1) = A,$ $\mathrm{supp}(f_2) = B,$ and $\mathrm{supp}(f_3) = C,$   we have  
$$
\prod_{e:\,e \in E} (f_1(e),f_2(e),f_3(e)) = \omega^{|A|}\omega^{|B|}\omega^{|C|} = \omega^{|A|+|B|+|C|}.  
$$
In addition, we have 
the following types of association between pairs and triples,  where $a$ and $b$ are nonzero elements of $\aA$: 
\begin{eqnarray*} 
(0,0)    &\longleftrightarrow&  (0,0,0); \,\, (a,a,a),\, a \neq 0 
\\ 
(0,a)    &\longleftrightarrow&  (a,a,0), \,\,(0,0,-a); \,\,(a+c,a+c,c), \,c \neq 0, -a  
\\
(a,0)    &\longleftrightarrow&  (a,0,0), \,\, (0,-a,-a);\,\,(a+c,c,c), \, c \neq 0, -a
\\
(a,b)    &\longleftrightarrow&  (a+b,b,0), \,\, (a,0,-b), \,\,(0,-a,-a-b);\,\,(a+b+c,b+c,c),\, c \neq 0, -b, -a-b. 
\end{eqnarray*}
Under the assumption that $2 \neq 0,$ $a \neq -a,$ 
\begin{eqnarray*} 
\sum_{a:\, a \in \aA}  \gamma(a,a,a) &=& 1 + (q-1)\omega^3, 
\\
\sum_{
(a,b,c): \, a,b,c \in \aA,\, a-b=0,\, c-b = a, a \neq 0} \gamma(a,b,c) &=&  
\omega + \omega^2 + (q-2)\omega^3, 
\\
\sum_{(a,b,c): \, a,b,c \in \aA,\, a-b=a,\, c-b = 0, a\neq 0} \gamma(a,b,c) 
&=& \omega + \omega^2 + (q-2)\omega^3, 
\\
\sum_{a:\, a,b,c \in \aA,\,a-b=a,\,b-c=b, a \neq 0, b \neq 0}  \gamma(a,b,c) &=& 3\omega^2 + (q-3)\omega^3. 
\end{eqnarray*}
The theorem now follows from Lemma 2.4.  
\qed
 
\vskip 0.2in 
The right-hand side of the equation in Theorem 5.1 cannot be manipulated so that we can use Lemma 2.2. 
However, over $\mathrm{GF}(2),$  the types of association degenerate on the right-hand side and Lemma 2.2 can be applied.   

\vskip 0.2in \noindent 
{\bf Theorem 5.2.}  Let $\sigma$ be an integer greater than $2,$ $\theta=2\rho\pi / \sigma,$ where $\rho$ is relatively prime to $\sigma,$ $\omega = \ee^{\iota \theta},$  $G$ be a rank-$r$ binary matrix with column set $E,$ and $\aA = \mathrm{GF}(2).$  Then 
$$
\sum_{k=0}^{\sigma-1} \omega^k |\mathcal{T} (k,\sigma)|
=
\omega^{|E| - r} (1 + \omega)^{|E|} (1 - \omega)^{2r} R\left( M(G);
\frac {2}{ \cos \theta - 1},
2(\cos \theta - 1)   \right).
$$

\vskip 0.1in\noindent 
{\bf Proof.}  Over $\mathrm{GF}(2),$ we have the simplest association of pairs with triples:    
\begin{eqnarray*} 
(0,0)    &\longleftrightarrow&  (0,0,0), \,\, (1,1,1)    
\\
(1,0)    &\longleftrightarrow&  (1,0,0), \,\, (0,1,1) 
\\  
(0,1)    &\longleftrightarrow&  (0,0,1), \,\, (1,1,0) 
\\
(1,1)    &\longleftrightarrow&  (1,0,1), \,\, (0,1,0). 
\end{eqnarray*} 
Using the same weight function $\gamma$ as in the proof of Theorem 5.1, we have 
$$
\gamma(0,0,0) + \gamma(1,1,1) = 1 + \omega^3  
$$
and 
$$
\gamma(1,0,0) + \gamma(0,1,1)  = \gamma(0,0,1) + \gamma(1,1,0)= \gamma(1,0,1) + \gamma(0,1,0)   
= \omega + \omega^2. 
$$  
Hence, by Lemma 2.4,  
$$
\sum_{
\vec{h}: \, \vec{h} =(h_1,h_2), \, h_1, \, h_2 \,\mathrm{are}\,\,\mathrm{flows}} 
(1 + \omega^3)^{|\vec{h}^{-1}(0,0)|}
 (\omega + \omega^2)^{|\vec{h}^{-1}(1,0)|+ |\vec{h}^{-1}(0,1)| + |\vec{h}^{-1}(1,1)|},  
$$
where the sum ranges over pairs $(h_1,h_2)$ of flows over $\mathrm{GF}(2)$ or equivalently, over (single) flows $\vec{h}$ over $\mathrm{GF}(2) \times \mathrm{GF}(2).$  Using Lemma 2.2 with 
$q = |\mathrm{GF}(2) \times \mathrm{GF}(2)| = 4,$  
together with the calculation 
$$
x - 1 = \frac {1 + \omega^3}{\omega + \omega^2} - 1  = \frac {1 - \omega + \omega^2}{\omega} - 1 =  \frac {1 - 2\omega + \omega^2}{\omega}  =
\omega^{-1} - 2 + \omega = 2 \cos \theta - 2,
$$
we obtain the equation in Theorem 5.2.    
\qed   

\vskip 0.2in\noindent 

We consider next special cases of Theorem 5.2 for small values of $\sigma.$  
Noting that over $\mathrm{GF}(2),$ 
$ \mathcal{T}(k,\sigma) $ is essentially the parcel 
$$\{(B_1,B_2,B_3): \, B_1, B_2, B_3\subseteq E, \, B_1 \Delta B_2, B_2\Delta B_3 \,\text{are unions of circuits},\, 
|B_1| +|B_2| +|B_3| \equiv k \,\text{mod}\,\sigma\},
$$
these cases are analogs of Theorems 16, 19, and 20 in \cite{Goodall}.   
When $\sigma = 3,$  $\omega^3 = 1,$ $(1 - \omega)^2(1 + \omega) =3,$ and $\omega(1 + \omega)= -1.$   
Theorem 5.2 yields the following result.   

\vskip 0.4in \noindent 
{\bf Corollary 5.3.} Let $G$ be a rank-$r$ binary matrix with column set $E$  and 
$\aA = \mathrm{GF}(2).$  Then  
$$
|\mathcal{T} (0,3)| - \tfrac {1}{2} |\mathcal{T} (1,3)| - \tfrac {1}{2} |\mathcal{T} (2,3)|
=
(-1)^{|E|-r}3^r R( M(G);-\tfrac {4}{3}, -3).
$$

\vskip 0.4in \noindent 
The next result follows from setting $\sigma = 4$  (so that $\omega = \iota,$ $(1-\omega)^2(1+\omega) = 2(1-\iota),$ and  $\omega(1 + \omega)= -(1-\iota)$)  in Theorem 5.2. 

\vskip 0.4in \noindent 
{\bf Corollary 5.4.}  Let $G$ be a rank-$r$ binary matrix with column set $E$  and 
$\aA = \mathrm{GF}(2).$  Then 
$$
|\mathcal{T} (0,4)| + \iota |\mathcal{T} (1,4)| - |\mathcal{T} (2,4)| - \iota |\mathcal{T} (3,4)|
=
(-1)^{|E|-r}2^r (1-\iota)^{|E|} R( M(G);-2, -2).
$$

\vskip 0.4in \noindent 
When $\sigma = 6,$  Theorem 5.2 gives an interpretation for the evaluation of the characteristic polynomial of a binary matroid at $4.$        

\vskip 0.4in \noindent 
{\bf Corollary 5.5.}
Let $G$ be a rank-$r$ binary matrix with column set $E$  and 
$\aA = \mathrm{GF}(2).$   If $|E|$ is even, 
$$
|\mathcal{T} (0,6)| + |\mathcal{T} (1,6)| - |\mathcal{T} (3,6)| -|\mathcal{T} (4,6)| 
=
(-3)^{|E|/2} \chi (M(G);4)   
\eqno(5)$$
and    
$$
|\mathcal{T} (1,6)| + |\mathcal{T} (2,6)| - |\mathcal{T} (4,6)| - |\mathcal{T} (5,6)| = 0.      
\eqno(6)$$
If $|E|$ is odd,   
\begin{align*}  
\tfrac {1}{2}[|\mathcal{T} (1,6)| + |\mathcal{T} (2,6)| - |\mathcal{T} (4,6)| - |\mathcal{T} (5,6)|] 
&= |\mathcal{T} (0,6)| + |\mathcal{T} (1,6)| - |\mathcal{T} (3,6)| -|\mathcal{T} (4,6)| 
\\ 
&=  
(-3)^{(|E|-1)/2} \chi (M(G);4)   
\end{align*}  
and
$$
|\mathcal{T} (0,6)| + \tfrac{1}{2}|\mathcal{T} (1,6)| - \tfrac {1}{2}|\mathcal{T} (2,6)|- |\mathcal{T} (3,6)| -\tfrac {1}{2}|\mathcal{T} (4,6)|+\tfrac {1}{2}|\mathcal{T} (5,6)| = 0.
$$

\vskip 0.1in\noindent
{\bf Proof.}  If $\omega = e^{\pi \iota /3},$ then 
$(1 - \omega)^2(1 + \omega) =-\sqrt{3} \iota$ and $\omega(1 + \omega)= \sqrt{3}\iota.$  
Hence, by Theorem 5.2, 
$$
\sum_{k=0}^5  \ee^{k \pi \iota /3}|\mathcal{T} (k,6)|
=
(-1)^r (\sqrt{3} \iota)^{|E|} R( M(G);-4, -1) = (\sqrt{3} \iota)^{|E|}\chi (M(G);4).
$$
When $|E|$ is even, the right-hand side is real and we have two equations: 
\begin{eqnarray*} 
|\mathcal{T} (0,6)| + \tfrac{1}{2}|\mathcal{T} (1,6)| - \tfrac {1}{2}|\mathcal{T} (2,6)|- |\mathcal{T} (3,6)| -\tfrac {1}{2}|\mathcal{T} (4,6)|  +\tfrac {1}{2}|\mathcal{T} (5,6)|
&=&
(-3)^{|E|/2}\chi (M(G);4),
\\
|\mathcal{T} (1,6)| + |\mathcal{T} (2,6)| -|\mathcal{T} (4,6)|   
- |\mathcal{T} (5,6)|
&=& 0.
\end{eqnarray*} 
The second equation yields Eq.~(6).  Adding Eq.~(6) to the first equation yields Eq.~(5).  The proof when $|E|$ is odd is similar.    
\qed 

\vskip 0.2in
Corollary 5.5 gives a parcel-size condition involving $6$ parcels for the existence of a nowhere-zero $4$-flow in a binary matroid.  This condition is the initial case of a sequence of parcel-size conditions involving $2m+2$ parcels for the existence of a nowhere-zero  $2^m$-flow.       

\vskip 0.2in \noindent 
{\bf Theorem 5.6.}  Let $m$ and $\tau$ be integers, $m \geq 2,$ $\tau \ge 1,$ $\omega = \ee^{\pi \iota /\tau},$ $G$ be a rank-$r$ binary matrix, and 
$$
\mathcal{T}_m (k,2\tau)  =  \{ (f_1,f_2,\ldots,f_{m+1}):\, f_j:E \to \mathrm{GF}(2), \, f_j - f_{j+1}  \,\,\text{are flows}, \,\,
\sum_{j=1}^{m+1} |\mathrm{supp}(f_j)| \equiv k \,\,\mathrm{mod}\, 2 \tau \}.
$$
Then (when $\tau = m+1$),  
$$ 
\sum_{k=0}^{2m+1} \ee^{k \pi \iota/(m+1)} |\mathcal{T}_m (k,2m+2)| \neq 0
\,\,\,\,\text{if and only if} \,\,\,\,\chi (M(G); 2^m) \neq 0. 
$$  

\vskip 0.1in\noindent 
{\bf Proof.}
 Let $(a_1,a_2, \ldots, a_m)$ be an $m$-tuple in $\mathrm{GF}(2)^m.$  Then there are two $(m+1)$-tuples associated with it:  $\vec{b}$ and its complement $\vec{b} +\vec{1},$ where 
$$\vec{b} = (a_1+a_2+ \cdots + a_{m-1} + a_m, a_2+a_3+ \cdots + a_{m-1}+ a_m, \ldots, a_{m-1} +a_m,a_m,0) 
\eqno(7)$$ 
and $\vec{1}$ is the vector with all coordinates equal to $1.$  Let $\wt(\vec {b})$ be the number of $1$'s in $\vec{b}.$  

Let $\omega =  \ee^{\pi \iota /\tau}$ and choose the weight function 
$$
\gamma(\vec{b}) =  \omega^{\wt(\vec{b})}.    
$$
By  Lemma 2.4,    
$$
\sum_{k=0}^{2 \tau-1} \omega^k |\mathcal{T}_m(k, 2 \tau)| 
=
\sum_{\vec{h}: \,\vec{h} \,\,\text{is a flow over} \,\,\mathrm{GF}(2)^m} 
\,\,
(1 + \omega^{m+1})^{|\vec{h}^{-1}(\vec{0})|} \prod_{\vec{a}:\vec{a} \neq \vec{0}}  
(\omega^{\wt (\vec{b})} + \omega^{m + 1 - \wt (\vec{b})})^{|\vec {h}^{-1}(\vec {a})|},
$$
where for a nonzero $m$-tuple $\vec{a}$ in the product on the right-hand side, $\vec{b}$ is the $(m+1)$-tuple defined in Eq.~(7).

Now let $\tau = m+1$ and $\omega = \ee^{\pi \iota/(m+1)}.$   Then $1 + \omega^{m+1} = 0$ and 
$$ 
\omega^{\wt (\vec {a})} + \omega^{m + 1- \wt (\vec {a})} = 
\omega^{\wt (\vec {a})} - \omega^{- \wt (\vec {a})} = 
2 \iota \sin \left(\frac {\wt (\vec{a}) \pi}{m+1} \right) . 
$$
Hence, 
\begin{eqnarray*}
\sum_{k=0}^{2m+1} \ee^{k \pi \iota/(m+1)}| \mathcal{T}_m(k,2m+2) | 
&=& 
\sum_{\vec{h}: \,\vec{h} \,\,\text{is a flow over} \,\,\mathrm{GF}(2)^m}
\,\, 
0^{|\vec{h}^{-1}(\vec{0})|} \, \prod_{\vec{a}:\, \vec{a} \neq \vec{0}}  
\left(  2 \iota \sin \left(\frac {\wt (\vec{a}) \pi}{m+1} \right)  \right)^{|\vec{h}^{-1}(\vec {a})|} 
\\
&=& 
(2 \iota)^{|E|} \sum_{\vec{h}: \,\vec{h} \,\,\text{is a flow over} \,\,\mathrm{GF}(2)^m,\,
\vec{h}^{-1}(\vec{0}) = \emptyset}
\quad \prod_{\vec{a}:\, \vec{a} \neq \vec{0}}  
\left( \sin \left(\frac {\wt (\vec{a}) \pi}{m+1} \right)   \right)^{|\vec{h}^{-1}(\vec {a})|}.
 \end{eqnarray*}
Since 
$$
0 < \sin \left(\frac {\wt (\vec{a}) \pi}{m+1} \right) \leq 1
$$ 
when $\vec{a} \neq \vec{0},$  the sum is over positive terms, one for each nowhere-zero flow  $\vec{h}$ 
over $\mathrm{GF}(2)^m.$  Hence,  
the sum is nonzero if and only if there is a nowhere-zero flow over $\mathrm{GF}(2)^m,$ that is to say, 
if and only if $\chi (M(G); 2^m) \neq 0.$  
\qed

\vskip 0.4in\noindent
{\bf\large  6.  Parcels defined using inner products}

\vskip 0.2in\noindent
Let $f, g:E \to \mathrm{GF}(q).$   Their {\sl inner} or {\sl dot product} $\langle f,g \rangle$ is defined by
$$
\langle f,g \rangle = \sum_{e:\, e \in E} f(e)g(e).
$$ 
In this section, we consider parcels in $\mathcal{F}$ defined by inner products.  For $G$ a matrix over $\mathrm{GF}(q)$  and $a$ an element in $\mathrm{GF}(q),$
let
\begin{eqnarray*} 
\mathcal{M}^{\bullet}(a, q)
=  
\{ (f,g): f - g \,\,\mathrm{is} \,\,\mathrm{a} \,\,\mathrm{flow}, \,
\langle f,g \rangle = a\}.  
\end{eqnarray*}
We shall restrict our attention to fields of prime order.  

\vskip 0.2in\noindent 
{\bf Theorem 6.1.}  Let $p$ be an odd prime, $\omega$ be a primitive $p$-th root of unity, and $G$ be a matrix over $\mathrm{GF}(p),$ where $p$ is an odd prime.  Then 
$$
\sum_{a:\,a \in \mathrm{GF}(p)} \omega^a |\mathcal{M}^{\bullet}(a,p)| 
=
\Omega^{|E|} \sum_{h:\, h\,\,\mathrm{is} \,\,\mathrm{a} \,\,\mathrm{flow}} 
\omega^{-\frac {1}{4} \langle h,h \rangle},  
$$ 
where 
$$
\Omega = 1 + \sum_{b:\,b\,\,\mathrm{is} \,\,\mathrm{a} \,\,\mathrm{square}} 2\omega^b  
$$
and the exponent $-\frac {1}{4} \langle h,h \rangle$ is evaluated in $\mathrm{GF}(p)$ and interpreted as an 
integer modulo $p.$  
\vskip 0.1in\noindent 
{\bf Proof.}  We use the weight function $ \gamma (b,c) = \omega^{bc}. $  Then 
$$ 
\sum_{b,c:\, b-c = 0} \gamma(b,c) = \sum_{b=0}^{p-1}  \gamma(b,b) = \sum_{b=0}^{p-1}  \omega^{b^2} = \Omega, 
$$ 
and for $a \neq 0,$  
$$
\sum_{b,c:\, b-c = a} \gamma (b,c) = \sum_{b=0}^{p-1} \omega^{b(b-a)} 
= \omega^{-a^2/4} \sum_{b=0}^{p-1} \omega^{b^2 - ba + (a/2)^2}= 
\omega^{-a^2/4} \sum_{b=0}^{p-1} \omega^{(b - a/2)^2} 
= \omega^{-a^2/4} \Omega.
$$
Hence, by Lemma 2.4, 
$$
\sum_{a:\,a \in \mathrm{GF}(p)} \omega^a |\mathcal{M}^{\bullet}(a,p)| 
=
\Omega^{|E|} \sum_{h:\, h\,\,\mathrm{is} \,\,\mathrm{a} \,\,\mathrm{flow}} 
\,\,\,\,\prod_{a:\, a\in \mathrm{GF}(p)}(\omega^{-a^2/4})^{|h^{-1}(a)|}.  
$$ 
To finish the proof, we use the fact that 
$$
\sum_{a:\, a \in \mathrm{GF}(p)} a^2 |h^{-1}(a)| = \langle h,h \rangle.  
$$
\qed 

\vskip 0.2in\noindent 
We note that because $1+\omega+\omega^2+ \cdots+\omega^{p-1} = 0,$ we can write 
$$
\Omega = \sum_{a=1}^{p-1} \left( \frac {a}{p} \right) \omega^a,
$$ 
where the LeGendre symbol $(\frac {a}{p})$ equals $1$ if $a$ is a square and $-1$ if $a$ is not a square. In other words, $\Omega$ is a Gauss sum.  In particular, it is known that the absolute value of $\Omega$ equals $\sqrt{p}.$  (See, for example, \cite{IRosen}, Chap.~8.) 

We can apply Lemma 2.2 to Theorem 6.1 to obtain an evaluation of the rank generating polynomial 
for only one case, that of ternary matroids.

\vskip 0.2in\noindent
{\bf  Corollary 6.2.}
Let $G$ be a rank-$r$ ternary matrix on the column set $E.$     Then
$$ 
|\mathcal{M}^{\bullet}(0,3)| + \ee^{2 \pi \iota /3} |\mathcal{M}^{\bullet}(1,3)| + 
\ee^{-2 \pi \iota /3} |\mathcal{M}^{\bullet}(-1,3)|
= 
\ee^{(5r - |E|) \pi \iota /6}   \sqrt{3}^{|E|+r}
R \left( M(G);  \sqrt{3} \ee^{ -5\pi \iota/6}, \sqrt{3} \ee^{5 \pi \iota /6}\right).
$$

\vskip 0.1in\noindent
{\bf Proof.} 
When $p=3,$  $\omega = \ee^{2 \pi \iota/3}$ and for 
$a \in \mathrm{GF}(3),$ 
$$
\omega^{-a^2/4} = \left\{\begin{array}{cc} 1 \,\,  &\text{if} \,\,a = 0,  \\ 
\omega^{-1} \,\, &\text{if}\,\, a \neq 0. 
\end{array} \right.  
$$  
By Theorem 6.1 and Lemma 2.2,   
\begin{eqnarray*}
|\mathcal{M}^{\bullet}(0,3)| + \omega |\mathcal{M}^{\bullet}(1,3)| + 
\omega^{-1}  |\mathcal{M}^{\bullet}(-1,3)| 
& =& 
\Omega^{|E|} \sum_{h:\, h\,\,\text{is  a flow}} 
\,\,\,\,  1^{|h^{-1}(0)|} \omega^{- [ |h^{-1}(1)| + |h^{-1}(-1)| ]}  
\\
&= &
\Omega^{|E|} \omega^{-|E|} \sum_{h:\, h\,\,\text{is a flow}} 
\,\,\,\,  \omega^{|h^{-1}(0)|} 
\\ 
&= &
\Omega^{|E|} \omega^{-|E|}(\omega - 1)^r  R\left(M(G); \frac {3}{\omega - 1}, \omega - 1\right).  
\end{eqnarray*} 
To finish the proof,  use $\omega = \ee^{2 \pi \iota/3},$  $\Omega = 1 + 2 \omega = \sqrt{3} \iota,$  and $\omega-1 = \sqrt{3} \ee^{5 \pi \iota /6}.$        
\qed

\vskip 0.2in\noindent
The value $R(M(G); \sqrt{3} \ee^{-5\pi \iota/6}, \sqrt{3} \ee^{5 \pi \iota /6})$ (which is the value 
$T(M(G);\ee^{-2\pi\iota/3}, \ee^{2 \pi \iota/3})$ of the Tutte polynomial)  was calculated, up to a root of unity by Jaeger \cite{Jg} and exactly by Gioan and Las Vergnas \cite{GLaV}.  It is a root-of-unity multiple of $\sqrt{3}^d,$ where $d$ is the bicycle dimension of $G.$  (The {\sl bicycle dimension} is the dimension of the vector space 
$U \cap V,$ where $U$ is the row space of $G$ and $V$ is the row space of a matrix $H$ orthogonally dual to $G.$ )  This is shown by evaluating the sum 
$\sum_{h} \omega^{\langle h,h \rangle}$ for a ternary matrix $G.$  With the knowledge that $|\Omega| = \sqrt{p},$ one can easily extend the elegant algebraic argument of Gioan and Las Vergnas to obtain the following result: if $p$ is an odd prime, $\omega$ is a primitive $p$-th root of unity,  and $G$ is a rank-$r$ $\mathrm{GF}(p)$-matrix, then 
$$
\left|\sum_{h:\, h\,\,\mathrm{is} \,\,\mathrm{a} \,\,\mathrm{flow}}  \omega^{\langle h,h \rangle} \right| 
= \sqrt{p}^{r+d}.
$$  
However, when $p > 3,$   our method fails to convert this sum to an evaluation of the rank generating polynomial.        

Theorem  6.1 covers the cases of odd primes. For completeness, we describe what happens when $p=2.$   Over $\mathrm{GF}(2),$ $\langle f,g \rangle \equiv |\mathrm{supp}(f) \cap \mathrm{supp}(g)|\,\,\mathrm{mod} \,\,2.$  Hence, $|\mathcal{M}^{\bullet}(k,2)| = |\mathcal{M}^{\cap}(k,2)|.$   By Corollary 4.17, if  $G$ is a binary matrix, then   
$$ 
|\mathcal{M}^{\bullet}(0,2)| - |\mathcal{M}^{\bullet}(1,2)| 
= 
(-1)^r  2^{|E|}  R ( M(G);  -2, - 1) = 2^{|E|}\chi (M(G);2).
$$   

\vskip 0.4in\noindent
{\bf\large   7.  Parcel-weight enumerators}

\vskip 0.2in\noindent
To end this paper, we consider parcels defined by ``letting $\sigma$ go to infinity''.  
When $\sigma$ is sufficiently large (compared with $|E|,$ $r,$ and $q$),  then the congruence condition becomes an absolute condition, in the sense that   
``is congruent to $k$ modulo $\sigma$'' strengthens to ``is equal to $k$''.   As a typical example, we consider the parcels  in Theorem 4.7.  Define $\mathcal{M}_{1,-1}(k,\infty)$ by     
$$  
\mathcal{M}_{1,-1}(k,\infty) = \{ (f,g): \, (f,g) \in \mathcal{F},  \,\,
|\mathrm{supp}(f)| - |\mathrm{supp}(g)| = k\}.
$$
Since Lemmas 2.2 and 2.3 are formal, we can replace the root $\omega$ of unity by an indeterminate $X$ in Theorem 4.7 and its proof.  This allows us to express the {\sl (parcel-weight) enumerator}  $\sum_{k} |\mathcal{M}_{1,-1}(k,\infty)| X^k$ 
as a specialization of the rank generating polynomial.    
In general, the weight enumerator is a Laurent polynomial, that is, a polynomial in $X$ and $X^{-1}.$  

\vskip 0.2in\noindent
{\bf  Theorem 7.1.}   Let $X$ be an indeterminate, $G$ a rank-$r$ totally-$\mathbb{P}$ matrix with column set $E,$  
and $\aA$ a $\mathbb{P}$-module of order $q.$  Then 
\begin{eqnarray*} 
&& \quad \sum_{k:\, \, k \,\,\mathrm{is} \,\,\mathrm{an} \,\,\mathrm{integer}} |\mathcal{M}_{1,-1}(k,\infty)| X^k 
\\
&&= 
(X - 2 + X^{-1})^r (X - 2 + X^{-1} + q)^{|E|-r}
R \left( M(G);  \frac {q(X - 2 + X^{-1} + q)}
{X - 2 + X^{-1}}, 
\frac {X - 2 + X^{-1}}
{X - 2 + X^{-1} + q} \right).
\end{eqnarray*}

\vskip 0.3in\noindent 
{\bf Corollary 7.2.}  The sizes $|\mathcal{M}_{1,-1}(k,\infty)|,$ and hence, the sizes 
$|\mathcal{M}_{1,-1}(k,\sigma)|,$ depend only on the matroid $M(G).$ 
 
\vskip 0.2in 
Most of the theorems in this paper have enumerator versions. The enumerator version of Theorem 3.1 extends a theorem of Greene (Corollary 4.5, \cite{Greene}; see \cite{Farr} for background and related results).  To see this, let 
$$
\mathcal{H}(k,\sigma) = \{h:\, h \,\,\mathrm{is} \,\, \mathrm{a} \,\, \mathrm{flow},\,
|\mathrm{supp}(h)| \equiv k\,\,\mathrm{mod}\, \sigma\}.  
$$
The sets $\mathcal{H}(k,\sigma)$ are not parcels, but $ |\mathcal{O}(k,\sigma)| = q^{|E|} |\mathcal{H}(k,\sigma)| $ 
(see Section 3).  
Hence, we have the following theorem.   

\vskip 0.2in\noindent
{\bf  Theorem 7.3.}  Let $\sigma$ be an integer greater than $1,$ $\omega$ be a primitive $\sigma$-root of unity, $X$ be an indeterminate, $G$ be a rank-$r$ totally-$\mathbb{P}$ matrix with column set $E,$  
and $\aA$ be a $\mathbb{P}$-module of order $q.$   Then 
\begin{eqnarray*}
\sum_{k=0}^{\sigma - 1}  \omega^k |\mathcal{H}(k,\sigma)|
&=&
\omega^{|E|-r} (1 - \omega)^r R \left( M(G); \frac {q \omega}{1 -\omega} , \frac {1-\omega}{\omega} \right) ,
\\
\sum_{k:\, k \ge 0}  |\mathcal{H}(k,\infty)| X^k 
&=&
X^{|E|-r} (1 - X)^r R \left( M(G); \frac {q X}{1 -X} , \frac {1-X}{X} \right).  
\end{eqnarray*}

\vskip 0.4in\noindent 
When $G$ is a matrix over $\mathrm{GF}(q)$ and $\aA = \mathrm{GF}(q),$ the homogeneous bivariate polynomial 
$\sum_{k:\, k \ge 0}  |\mathcal{H}(k,\infty)| X^k Y^{|E| - k}$ is the weight enumerator of the linear code 
generated by the matrix $G$ as defined, say, in \cite{MacW}. Thus, in this case, the second part of Theorem 7.3 is Greene's theorem.     

\vskip 0.4in\noindent
{\bf Acknowledgment.}  We thank the referees for many comments and suggestions.

\end{document}